\documentclass[10pt]{article}
\usepackage{amsmath}
\usepackage{amsfonts}
\usepackage{amssymb}
\usepackage{graphicx}

\newtheorem{definition}{Definition}
\newtheorem{theorem}{Theorem}
\newtheorem{corol}{Corollary}
\newtheorem{lemma}{Lemma}
\newtheorem{remark}{Remark}
\begin{document}
\title{An Optimal  Control Problem for the Steady Nonhomogeneous
Asymmetric Fluids}
\author{Exequiel Mallea-Zepeda$^1$, Elva Ortega-Torres$^2$, \'Elder J. Villamizar-Roa$^3$}
\date{\small$^1$\it Departamento de Matem\'atica, Universidad de Tarapac\'a, Arica, Chile\\
$^{2}$\it Departamento de Matem\'aticas, Universidad Cat\'olica
del Norte, Antofagasta, Chile\\
$^3$\it Escuela de Matem\'aticas, Universidad Industrial de Santander, Bucaramanga, Colombia}
\maketitle

\footnotetext{$^1$  E-mail:{\tt emallea@uta.cl}}
\footnotetext{$^2$E-mail: {\tt eortega@ucn.cl}}
\footnotetext{$^3$E-mail: {\tt jvillami@uis.edu.co}}
\date{}

\begin{abstract}
We study an optimal boundary control problem for the two-dimensional stationary micropolar fluids system with variable density. We control the system by considering boundary controls, for the velocity vector and angular velocity of rotation of particles, on parts of the boundary of the flow domain. On the remaining part of the boundary, we consider mixed boundary conditions for the vector velocity (Dirichlet and Navier conditions) and Dirichlet boundary conditions for the angular velocity. We analyze the existence of a weak solution obtaining the fluid density as a scalar function of the stream function. We prove the existence of an optimal solution and, by using the Lagrange multipliers theorem, we state first-order optimality conditions.  We also derive, through a penalty method, some optimality conditions satisfied by the optimal controls.\newline

\end{abstract}

\section{Introduction}
Incompressible fluids with variable density (non homogeneous Navier-Stokes equation) correspond to a coupling between the equation for the velocity
given by the conservation of momentum, the transport equation for the density provided by the mass conservation law,
and the incompressibility condition. This kind of fluids are relevant to be analyzed from the mathematical and physical point of view. They can be used to model, among others, stratified fluids \cite{Pedlosky}, meeting of fluids coming from various regions
with different densities, like the  junction of pipes filled with incompressible fluids with different densities or the
junction of two or more rivers \cite{Santos1}. There exists a considerable number of papers devoted to the mathematical
analysis of the non homogeneous Navier-Stokes equations, principally in the non stationary case, including results
when the initial density is assumed to be positive or when the initial-vacuum is allowed (see \cite{LGE, Lions, Simon}
and references therein); however, not much is known about the stationary case including optimal control problems, where the state equations are given by the equations describing the motion of a viscous incompressible fluid with variable density.

An important model which generalizes the non homogeneous Navier-Stokes equation is given by the non homogeneous micropolar fluids. Non homogenous micropolar fluids refer to the micropolar fluid model with variable density; meanwhile, micropolar fluids are fluids with microstructure and asymmetric stress tensor. Physically, they represent fluids consisting of
randomly oriented (or spherical)  particles suspended in a viscous medium, when the
deformation of fluid particles is ignored \cite{Eringen1}. This model, in stationary state,
is given by the following system of partial differential equations which expresses  the balance of momentum, mass, and moment of momentum
(cf. \cite{Eringen2,lukaszewicz}):
\begin{equation}
\left\{
\begin{array}{rcl}
-(\mu+\mu_r)\Delta{\bf u}+ \rho({\bf u} \cdot\nabla) {\bf u} +\nabla p &=&
2\mu_r{\rm rot}\, {\bf w}+\rho{\bf f}\ \mbox{in}\ \Omega,\vspace{0.1cm}\\
-(c_a+c_d)\Delta {\bf w}-(c_0+c_d-c_a)\nabla{\rm div}\,{\bf w}+\rho({\bf u}\cdot\nabla ){\bf w}
+4\mu_r{\bf w}&=&2\mu_r{\rm rot}\, {\bf u}+\rho{\bf g} \ \mbox{in}\  \Omega,\vspace{0.1cm}\\
{\bf u}\cdot\nabla\rho&=&0\ \mbox{ in }\ \Omega,\vspace{0.1cm}\\
\rm{div}{\bf u} &=&0\ \mbox{in}\ \Omega,\\
\end{array}
\right.  \label{eq1}
\end{equation}
where $\Omega$ is a  connected bounded domain of
$\mathbb{R}^3$ with Lipschitz boundary, ${\bf u}=(u_1,u_2,u_3)$ is the velocity field, $\rho$ denotes the density, $p$
represents the pressure, and ${\bf w}$ is the microrotation field interpreted as the angular velocity field of rotation
of particles. The fields ${\bf f}=({\rm f}_1,{\rm f}_2,{\rm f}_3)$ and ${\bf g}=({\rm g}_1,{\rm g}_2,{\rm g}_3)$ represent
external sources of linear and
angular momentum respectively. The positive constants
$\mu,\mu_r,c_a,c_d,c_0$ characterize isotropic properties
of the fluid; in particular, $\mu$ denotes the dynamic viscosity, and $\mu_r,c_a,c_d, c_0$ are new viscosities
connected with physical characteristics of the fluid. These constants satisfy $c_0+c_d>c_a$. For simplicity we
denote $\mu_1=\mu+\mu_r$, $\mu_2=c_a+c_d$ and $\mu_3=c_0+c_d-c_a$. When the microrotation viscous effects
are neglected, that is $\mu_r=0$, or the microrotation velocity is null, the micropolar fluid
model reduces to the classical incompressible Navier-Stokes system.

From the mathematical point of view, the micropolar
fluid system, with constant density, has been studied by several authors, and important results
on well-posedness, large time asymptotic behavior and general
qualitative analysis, have been obtained (see, for instance,
\cite{FerreVilla,lukaszewicz,VillaRod,VillaRodRojas}
and references therein). However, as far as we know, the variable density stationary model (\ref{eq1})  has only been previously considered  in \cite{Vitoriano},
where, by using the Galerkin method, the author proved the existence of weak solutions
for the system (\ref{eq1.1})-(\ref{eq1.3}) with Dirichlet boundary conditions. The main difficulty of studying model (\ref{eq1}) is due to the first-order equation ${\bf u}\cdot\nabla \rho=0$ in $\Omega$
with $\rho=\rho_g$ on $\Gamma_0\subseteq\partial\Omega.$ Even in the particular case when ${\bf w}={\bf 0},$ there are fewer results available in the literature related to
the existence of solutions for (\ref{eq1}), and they depend on the
dimension of the domain $\Omega$ (see \cite{Santos2, Santos3, Santos4, Santos5, Frolov1, Frolov2,
Santos1}). In particular, in
\cite{Frolov1} the author proved the existence of a solution
$[{\bf u},\rho]$ for system (\ref{eq1}) in the class ${\bf H}^2(\Omega)\times
C^\theta(\Gamma_0)$ provided that ${\bf f}\in {\bf L}^2(\Omega),$ ${\bf u}_{\bf h}\in
{\bf C}^2(\Gamma),$ and $\rho_g\in C^\theta(\Gamma_0)$ for $\theta>0,$ with ${\bf u}_{\bf h}$ and $\rho_g$ being the boundary data for the velocity and density respectively. This result was
improved in \cite{Santos1} where the existence of a weak solution
with boundary values for the density prescribed in $L^\infty$ was
obtained. Still in the case 2D, but in unbounded domains,
some results related to the Leray problem have been obtained in \cite{Santos2, Santos3,
Santos4, Santos5}. The existence of solutions in the case 3D seems to be
more difficult to handle and, differently to the non stationary case,
we only know the paper \cite{Frolov2}.

In this paper, we confine ourselves to two-dimensional flows in a bounded domain with boundary $\Gamma$ of class $C^2$. Such a flow can be interpreted as being a
cross section of the three-dimensional domain $\Omega$ by making $x_3=c$, where $c$ is a constant. In this case,
it is assumed that the velocity component $u_3$ in the $x_3-$direction is zero, and the axes of rotation of particles
are parallel to the $x_3-$axis. Then, for ${\bf x}=[x_1,x_2]\in\Omega\subset\mathbb{R}^2$, the fields ${\bf u}, {\bf w},
\rho$ and $ p$ reduce to ${\bf u}=[u_1({\bf x}),u_2({\bf x}),0]$, ${\bf w}=[0,0, {\rm w}_3({\bf x})]$,
$\rho=\rho({\bf x})$, and $p=p({\bf x})$. Also, the external sources can be written as ${\bf f}=[f_1({\bf x}), f_2({\bf x}),0]$
and ${\bf g}=[0,0,{\rm g}_3({\bf x})]$. Consequently, from now on we assume the following notations: ${\bf u}=[u_1({\bf x}),u_2({\bf x})]$,
${\rm w}={\rm w}_3({\bf x})$, $\rho=\rho({\bf x})$, $p=p({\bf x})$, ${\bf f}=[f_1({\bf x}),f_2({\bf x})]$, and
${\rm g}={\rm g}_3({\bf x})$. Then, by observing that
\[
{\rm rot}\, {\bf u}=\frac{\partial u_2}{\partial x_1}-\frac{\partial u_1}{\partial x_2}, \quad
{\rm rot}\, {\rm w}=[\frac{\partial {\rm w}}{\partial x_2}, -\frac{\partial {\rm w}}{\partial x_1}],
\quad {\rm div}\, {\rm w}=0,
\]
and considering ${\bf u}$, ${\rm w}$, $\rho$, $p$, ${\bf f}$, and ${\rm g}$ in the system (\ref{eq1}),
we obtain the following two-dimensional system
\begin{equation}
\left\{
\begin{array}{rcl}
-\mu_1\Delta{\bf u}+ \rho({\bf u} \cdot\nabla) {\bf u} +\nabla p &=&
2\mu_r{\rm rot}\, {\rm w}+\rho{\bf f}\ \mbox{in}\ \Omega,\vspace{0.1cm}\\
-\mu_2\Delta {\rm w}+\rho({\bf u}\cdot\nabla ){\rm w}
+4\mu_r{\rm w}&=&2\mu_r{\rm rot}\, {\bf u}+\rho{\rm g} \ \mbox{in}\  \Omega,\vspace{0.1cm}\\
{\bf u}\cdot\nabla\rho&=&0\ \mbox{ in }\ \Omega,\vspace{0.1cm}\\
\rm{div}\,{\bf u} &=&0\ \mbox{in}\ \Omega.\\
\end{array}
\right.  \label{eq1.1}
\end{equation}

In this paper, we prove the existence of weak solutions for (\ref{eq1.1}) and then, we study an
optimal boundary control problem where the state equations are given by the
weak solutions of (\ref{eq1.1}). For this purpuse
we consider the following boundary conditions:
\begin{equation}\label{eq2.1}
{\bf u}={\bf u}_{\boldsymbol{g}_1}:=\left\{
\begin{array}{ccl}
{\bf u}_0 & \mbox{on} & \Gamma_0,\\
\mbox{\bf\em g}_1 & \mbox{on} & \Gamma_1,\\
\end{array}
\right.\ \
\left[D({\bf u}){\bf n}+\alpha{\bf u}\right]_{\rm tang}= 0\ \mbox{on}\ \Gamma_2,\ \ {\bf u}\cdot {\bf n}=0\ \mbox{on}\ \Gamma_2,
\end{equation}
\begin{equation}\label{eq2.1b}
\int_\Gamma{\bf u}\cdot{\bf n}\,d\Gamma=0,\ \ \ \rho=\rho_0>0\ \mbox{ on }\Gamma_0,\ \ \
{\rm w}_{g_2}=\left\{
\begin{array}{ccl}
{\rm w}_0 & \mbox{on} & \Gamma_0,\\
\mbox{\em g}_2 & \mbox{on} & \Gamma_3.\\
\end{array}
\right.
\end{equation}
Here the boundary $\Gamma$ of $\Omega$ is of class $C^2$ and $\Gamma=\Gamma_0\cup\Gamma_1\cup \Gamma_2=\Gamma_0\cup\Gamma_3$,
where $\Gamma_0\cap\Gamma_i=\emptyset$, $i=1,3,$
$\overline{\Gamma}_1\cap\overline{\Gamma}_2=\overline{\Gamma}_0\cap\overline{\Gamma}_2=\emptyset.$  We assume  that
\begin{equation}\label{eq1.2}
\Gamma_0\mbox{ is an arcwise connected closed set on } \Gamma, \mbox{with measure} (\Gamma_0)>0.
\end{equation}
The parts $\Gamma_1$ and $\Gamma_3$ have nonempty interior, but $\Gamma_2$ may be an empty set.
The functions $\rho_0,{\bf u}_0$ and ${\rm w}_0$ are defined on $\Gamma_0,$ and
the functions $\mbox{\bf\em g}_1$, $\mbox{\em g}_2$ describe the Dirichlet
boundary control for ${\bf u}$ on $\Gamma_1$ and for ${\rm w}$ on
$\Gamma_3$ respectively. The controls $\mbox{\bf\em g}_1,\mbox{\em g}_2$ lie in closed convex sets
$\mathcal{U}_1\subset{\bf H}^{1/2}(\Gamma_1)$ and
$\mathcal{U}_2\subset{ H}^{1/2}(\Gamma_3)$ respectively. We assume that

\begin{equation}\label{eq1.3}
{\bf u}_0\cdot{\bf n}>0 \mbox{ on } \Gamma_0\ \mbox{(outflow)}\ \mbox{ or }\ {\bf u}_0\cdot{\bf n}<0\mbox{ on }\Gamma_0\ \mbox{(inflow)},
\end{equation}
where ${\bf n}$ denotes the outward normal vector on $\Gamma$.\newline
\\
\vspace{0.1cm}
\begin{minipage}{\linewidth}
\begin{center}
{\includegraphics[scale=0.8]{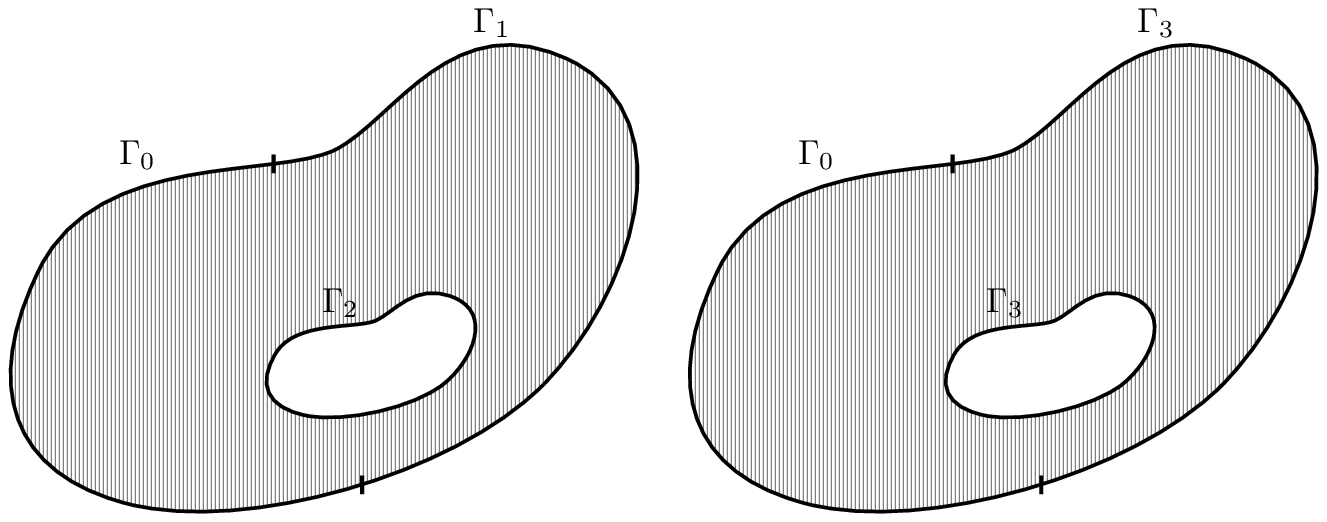}}\\
\end{center}
\end{minipage}
\begin{center}
\vspace{0.1cm}
\small{{\bf Fig. 1}{ Sketch of the domain $\Omega$} }
 \end{center}
The condition
$[D{(\bf u )}{\bf n}+\alpha {\bf u}]_{\rm tang}={\bf 0}, {\bf u}\cdot{\bf n} =0$ on $\Gamma_{2}$, is a Navier friction boundary condition. The term $[D({\bf u}){\bf n}+\alpha{\bf u}]_{\rm tang}:=
D({\bf u}){\bf n}+\alpha{\bf u}-[(D({\bf u}){\bf n}+\alpha{\bf u})\cdot{\bf n}]{\bf n}$ represents the tangential component of
the vector $D({\bf u}){\bf n}+\alpha{\bf u}$, where $D({\bf u}):=\frac12(\nabla{\bf u}+\nabla^T{\bf u})$ is the deformation
tensor, and $\alpha\geq0$ is the friction coefficient which measures the tendency of the fluid to slip
on $\Gamma_2$.
 The Navier boundary condition was proposed by Navier \cite{Navier},
who claimed that the tangential component of the viscous stress at the boundary should be proportional to
the tangential velocity. Navier boundary condition was also derived by Maxwell \cite{Maxwell} from the kinetic theory
of gases and rigorously justified as a homogenization of the no-slip condition on a rough boundary \cite{Jagger}.

We consider an objective functional  given by
a sum of functionals which measure, in the Lebesgue norm, the difference
between the velocity vector (respectively, the density and the microrotation velocity)
and a given prescribed velocity (respectively, a prescribed density and a
microrotation velocity). The objective functional also measures the
turbulence in the flow through a norm of the vorticity; it permits to describe the
resistance in the fluid due the viscous friction (see the cost functional in
(\ref{objective})). The state equations are given by  a weak formulation of the stationary micropolar fluids equations (\ref{eq1.1}) with boundary conditions (\ref{eq2.1})-(\ref{eq2.1b}). The exact mathematical
formulation will be given in Section 2.2 (see Definition \ref{weak}). The novelty of this paper lies in the following two aspects:
\begin{enumerate}
\item First,  we prove the existence of a weak solution for the stationary micropolar fluids equations
(\ref{eq1.1}) with boundary conditions (\ref{eq2.1})-(\ref{eq2.1b}). The definition of weak solution is given in Definition \ref{weak} in Subsection \ref{2.2} and he existence of weak solutions is given in Theorems \ref{teo_exist} and \ref{teo_exist_nohomo} in Section \ref{sec3}. We look for weak solutions with $\rho$ in the form $\rho=\eta(\psi):=\eta\circ\psi,$ for a continuous and positive function $\eta:\mathbb{R}\rightarrow\mathbb{R}$ and $\psi$ the stream function associated to the velocity vector, being this the reason why we consider $\Omega\subset \mathbb{R}^2$.
\item Second, we prove the
solvability of the optimal control problem. The existence of an optimal solution is given in Theorem \ref{teorem3} in Section \ref{sec4}. Posteriorly, by using the theorem of Lagrange multipliers, we state first-order optimality conditions. The first optimality conditions are obtained in Theorem \ref{Lagrange} in Section \ref{sec5}. We also
derive an optimality system in Remarks \ref{opt1} and \ref{opt2}.
\end{enumerate}
In order to obtain the first-order optimality conditions, we will use a penalty method.
This is a non standard technique which has been used previously in \cite{Casas, Illarionov,
Lee} to derive optimality conditions for optimal
control problems where the relation control-state is multivalued. To carry out this procedure, we introduce a family of penalyzed problems
which approximates the initial control problem (see Theorem \ref{theo_conver}); then, we analyze
their optimality conditions (see Theorem \ref{theo_mul-1}), and finally, we pass to the limit in the
parameter of penalization in order to derive the optimality
conditions of the original problem.

As far as we know, unlike the Navier-Stokes case, few works on optimal control problems
for micropolar fluids (with constant density) are available in the literature
 \cite{Exequiel,stavre3,stavre2,stavre1}. In \cite{stavre3}, a control problem for non stationary fluids in a two dimensional domains was analyzed;
in that paper, a viscosity
coefficient $\lambda,$ which achieves a desired field of the
microrotation velocity, is determined. In \cite{stavre2}, the author
studied an optimal control problem associated with the motion of a
micropolar fluid, with applications in the control of the blood
pressure. In \cite{stavre1}, the author analyzed, in a two-dimensional
domain, the relation between the microrotation and the vorticity of the
fluid. Recently, in \cite{Exequiel}, was considered an optimal boundary control problem for  micropolar fluids (with constant density)
equations in $3D$ bounded domains. Thus, the results of this paper can be seen as a 2D version of the results of
\cite{Exequiel} in the case of  micropolar fluids with variable density.


The outline of this paper is as follows: In Section 2, we establish the definition of weak solution and the optimal control problem to be considered.
In Section 3, we prove the existence of weak solutions. In Section 4, we prove the existence of an optimal solutions. In
Section 5, we derive first-order optimality conditions and, by
using the Lagrange multipliers theorem, we derive an optimality
system.

\section{Statement of the Problem}
\subsection{Function Spaces}\label{sub2.1}
Throughout this paper we will use the Sobolev space $H^1(\Omega),$ and
$L^p(\Omega),$ $1\leq p\leq \infty,$ with norms $\Vert
\cdot\Vert_{H^1}$ and $\Vert\cdot \Vert_p$ respectively. In
particular, the norm and inner
product in $L^2(\Omega)$ will be represented by $\Vert \cdot\Vert$ and
$(\cdot,\cdot)$ respectively. The norm $L^p(\Gamma)$ will be denoted by $\|\cdot\|_{L^p(\Gamma)}$.
Corresponding Sobolev spaces of
vector valued functions will be denoted by ${\bf H}^1(\Omega),$
${\bf L}^2(\Omega),$ and so on. We  will use the Hilbert space
${\ H}^1_0(\Omega)=\{{ u}\in { H}^1(\Omega)\,:\, { u}={ 0}\ \mbox{ on }\Gamma\}$, with the inner product
$({\rm u},{\rm v})_{{ H}^1_0}=(\nabla{\rm u},\nabla{\rm v})$ and the norm $\|{\rm u}\|_{H^1_0}=\|\nabla{\rm u}\|$.
We also consider the following solenoidal Banach spaces
${{\bf H}}_{\sigma}=\{{\bf u}\in {\bf H}^1(\Omega):
\mbox{div}\, {\bf u}=0 \mbox{ and }{\bf u}\cdot{\bf n}= 0\mbox{ on }\Gamma_2\},$
endowed  with the usual norm of ${\bf H}^1(\Omega),$ and  the space
$\tilde{\bf H}_\sigma=\{{\bf u}\in {\bf H}_\sigma: {\bf
u}={\bf 0}\mbox{ on }\Gamma\setminus\Gamma_2\},$
which is a Hilbert space with the inner product
$
({\bf u},{\bf v})_{\tilde{\bf H}_\sigma}=(D({\bf u}),D({\bf v}))
$
 and the norm
 $
 \|{\bf u}\|_{\tilde{\bf H}_\sigma}:=\|D({\bf u})\|.
 $

If $X$ is a general Banach space, its topological
dual will be denoted by $X'$ and the duality product by $\langle
\cdot,\cdot\rangle_{X'}$ or simply by $\langle \cdot,\cdot\rangle,$
unless this leads to ambiguity. The space ${ H}^{-1}$ denotes the dual of ${ H}^1_0(\Omega);$
the space ${\bf H}'$ denotes the dual of ${\bf H}_\sigma,$ and
the space $\tilde{\bf H}'$ denotes the dual of $\tilde{\bf H}_\sigma$.

If  $\Gamma_k$ is a connected subset of $\Gamma,$ we consider the trace space
${\bf H}^{1/2}(\Gamma_k)=\{ {\bf u}_{|_{\Gamma_k}}:{\bf u}\in {\bf H}^1(\Omega)\}$ (the restriction of
the elements of ${\bf H}^1(\Omega)$ to $\Gamma_k$) and
\begin{eqnarray*}
{\bf H}^{1/2}_{00}(\Gamma_k)&=&\{{\bf v}\in {\bf L}^2(\Gamma_k):\ \mbox{there exists}\ \hat{{\bf v}}
\in {\bf H}^{1/2}(\Gamma),\ \hat{\bf v}_{|_{\Gamma\setminus\Gamma_k}}
={\bf 0},\  \hat{{\bf v}}_{|_{\Gamma_k}}={\bf v}\}.
\end{eqnarray*}
In the case of scalar functions, we also use the space ${H}^{1/2}_{00}(\Gamma_k)$ which is defined similarly.
It can be verified that ${\bf H}^{1/2}_{00}(\Gamma_k)$ is a closed subspace
of ${\bf H}^{1/2}(\Gamma_k);$ moreover, ${\bf H}^{1}_{0}(\Gamma_{k})
\hookrightarrow {\bf H}^{1/2}_{00}(\Gamma_{k}) \hookrightarrow {\bf L}^{2}(\Gamma_{k})$
(cf. \cite{Dautray,de los Reyes1}).  For the space ${\bf H}^{1/2}(\Gamma_k),$
${\bf H}^{-1/2}(\Gamma_k)$
denotes its dual and $\langle \cdot,\cdot\rangle_{\Gamma_k}$ represents its duality product. The letter C will
denote diverse positive constants which may change from line to line or even within a same line.

For each function ${\bf u}\in{\bf H}_\sigma$ there exists  a scalar function $\psi\in H^2(\Omega)$ (stream-function)
such that
\begin{equation}\label{eq3}
{\bf u}={\rm rot}\,\psi=[\frac{\partial\psi}{\partial x_2},-\frac{\partial \psi}{\partial x_1}] \ \mbox{in} \ \Omega.
\end{equation}
Let $N:{\bf H}_\sigma\rightarrow H^2(\Omega)$ the linear operator
assigning to each  vector field ${\bf u}\in{\bf H}_\sigma$  its stream-function $\psi=N{\bf u}$
satisfying (\ref{eq3}).
The assumption ${\rm rot}\,\psi={\bf u}_0$ on $\Gamma_0$ implies that
\begin{eqnarray}
\frac{\partial \psi}{\partial {\bf n}}={\bf u}_0\cdot \boldsymbol{\tau},\ \ \frac{\partial \psi}{\partial \boldsymbol\tau}=-({\bf u}_0\cdot {\bf n})\ \mbox{ in }\ \Gamma_0,
\end{eqnarray}
where $\boldsymbol\tau$ denotes the outward tangent vector on $\Gamma.$
Thus, the boundary values of $\psi$ can be obtained by integrating with respect to the arc length, that is,
\begin{equation}\label{eq3z}
\psi({\bf x})= \int_{\Gamma_0({\bf x}_0,{\bf x})}
{\bf u}_0\cdot{\bf n}\,d\Gamma\quad {\bf x}\in\Gamma_0,
\end{equation}
where ${\bf x}_0$ is the initial point of the curve $\Gamma_0$ and $\Gamma_0({\bf x}_0,{\bf x})$ is the
part of the curve $\Gamma_0$ lying between the points ${\bf x}_0$ and ${\bf x}$ (cf. \cite{Frolov1}).
Notice that since $\psi\in H^2(\Omega),$ then  $\psi\in{ H}^{3/2}(\Gamma_0)\subset C^0(\Gamma_0)$; moreover, by (\ref{eq1.2})-(\ref{eq1.3}),
$\psi$ is strictly monotone on $\Gamma_0$. Therefore, there exists $\psi^{-1}\in C^0(\psi(\Gamma_0)),$ with
$\psi^{-1}:\psi(\Gamma_0)\subset\mathbb{R}\rightarrow \Gamma_0.$ Thus, if we assume that

\begin{equation}\label{bound_density}
\rho_0\in C^0(\Gamma_0),
\end{equation}
we can define the continuous function
$\tilde{\eta}(z)=\rho_0(\psi^{-1}(z)),$ $z\in\psi(\Gamma_0)\subset \mathbb{R}.$  Since  $\rho_0({\bf x})>0$ for all
${\bf x}\in\Gamma_0,$ and $\Gamma_0$ is an arcwise connected closed set in $\Gamma$, we can extend $\tilde{\eta}$ to $\mathbb{R}$ as a strictly positive scalar function $\eta$ such that
\begin{equation}\label{eq4}
\eta\in \mbox{ C}^0(\mathbb{R}), \ \ \eta(y) >0 \ \ \forall y\in\mathbb{R},
\ \eta(y)=\rho_0(\psi^{-1}(y)), \ \  y\in \psi(\Gamma_0).
\end{equation}
Therefore, under the above considerations, following \cite{Illarionov}, we define the density $\rho:\Omega\rightarrow\mathbb{R}$ as being
$$
\rho({\bf x})=\eta(\psi)({\bf x})=\eta(N{\bf u})({\bf x}),\ {\bf x}\in\Omega.
$$
\begin{remark}
By construction of $\eta$ it holds that $\eta(N{\bf u})=\rho_0$ on $\Gamma_0.$ Moreover, for ${\bf u}\in{\bf H}_\sigma$ and $\eta\in C^1(\mathbb{R})$ it holds
\begin{eqnarray}\label{densi}
\int_\Omega ({\bf u}\cdot\nabla \psi)\eta(N {\bf u})\,dx=0,\ \forall\psi\in H_0^1(\Omega).
\end{eqnarray}
If $\eta\in C^0(\mathbb{R}),$ we can regularize $\eta$ and the relation (\ref{densi}) remains true.
\end{remark}
\begin{remark}\label{remark2}
The operator $N:{\bf H}_\sigma\rightarrow H^2(\Omega)$  is continuous. The proof can be found in \cite{Illarionov}, Lemma 2.1.
\end{remark}
%
%
\subsection{Definition of Weak Solution}\label{2.2}
We consider the following operators
\begin{equation*}
\left\{\begin{array}{lll}
B:{\bf H}_\sigma\times {\bf H}_\sigma \times {\bf H}_\sigma \rightarrow \tilde{\bf H}', \quad
F:{\bf H}_\sigma \rightarrow \tilde{\bf H}',\\
\tilde{B}: {\bf H}_\sigma\times {\bf H}_\sigma\times H^1(\Omega) \rightarrow H^{-1}, \quad
G: {\bf H}_\sigma \rightarrow { H}^{-1},
\end{array}\right.
\end{equation*}
defined by
\begin{equation}\label{oper1}
\left\{
\begin{array}{lcl}
B[{\bf u},{\bf v},{\bf e}]=\eta(N{\bf u}){\bf v}\cdot\nabla{\bf e}, && F({\bf u})=\eta(N{\bf u}){\bf f},\\
\tilde{B}[{\bf u},{\bf v}, {\rm w}]=\eta(N{\bf u}){\bf v}\cdot\nabla {\rm w}, && G({\bf u})=\eta(N{\bf u}){\rm g}.
\end{array}
\right.
\end{equation}
For $\eta\in C^1(\mathbb{R})$ it holds that ${\rm div}(\eta(N{\bf u}){\bf u})={\rm div}(\rho{\bf u})={\bf u}\cdot \nabla\rho=0$; then
we get
\begin{equation}\label{eq6}
\langle B[{\bf u},{\bf u},{\bf v}], {\bf v}\rangle=0 \ \ \forall {\bf u} \in {\bf H}_\sigma,\,  {\bf v} \in\tilde{\bf H}_\sigma, 
\qquad \langle\tilde{B}[{\bf u},{\bf u}, {\rm w}], {\rm w}\rangle=0 \ \ \forall {\bf u}\in {\bf H}_\sigma,\, {\rm w}\in H^1_0(\Omega).
\end{equation}
For $\eta\in C^0(\mathbb{R}),$ by regularizing the function $\eta,$ the properties (\ref{eq6}) remain true.
Taking into account the operators defined in (\ref{oper1}), the problem (\ref{eq1.1})-(\ref{eq2.1b}) can be written as
\begin{equation}\label{eq7}
\left\{\begin{array}{rcl}
-\mu_1\Delta{\bf u}+ B[{\bf u},{\bf u},{\bf u}] +\nabla p&=& 2\mu_r{\rm rot}\, {\rm w}+F({\bf u})\ \mbox{ in } \Omega,\\
-\mu_2\Delta {\rm w}+\tilde{B}[{\bf u},{\bf u}, {\rm w}]+ 4\mu_r {\rm w} &=& 2\mu_r{\rm rot}\,{\bf u}+G({\bf u})\ \mbox{ in } \Omega,\\
\eta(N{\bf u})&=&\rho_0\ \mbox{ on }\Gamma_0,\\
{\bf u} &=& {\bf u}_{\boldsymbol{g}_1}\ \mbox{ on }\Gamma\setminus\Gamma_2,\\
{\bf u}\cdot{\bf n}&=&0\ \mbox{ on }\Gamma_2,\\
\left[D({\bf u}){\bf n}+\alpha{\bf u}\right]_{\rm tang}&=& 0\ \mbox{on}\ \Gamma_2,\\
{\rm w} &=&{\rm w}_{g_2}\ \mbox{ on }\Gamma.
\end{array}\right.
\end{equation}
\begin{lemma}(\cite{dragos,verfurth})\label{lema_partes}
Le ${\bf u}\in {\bf H}^2(\Omega)$ be a divergent free vector fields verifying the Navier boundary condition
and ${\bf v}\in {\bf H}^1(\Omega)$ is a divergence free vector field tangent to the boundary. Then,
$$
-\int_\Omega\Delta{\bf u}\cdot{\bf v}dx=2\int_\Omega D({\bf u}):D({\bf v})dx-2\int_\Gamma[D({\bf u}){\bf n}]_{\rm tang}\cdot{\bf v}d\,\Gamma.
$$
\end{lemma}
Through integration by parts and Lemma \ref{lema_partes}, we establish the following
definition of weak solution for system (\ref{eq1.1})-(\ref{eq2.1b}).
\begin{definition}\label{weak}
Let ${\bf f}\in{\bf L}^2(\Omega),\, {\rm g}\in L^2(\Omega)$, ${\bf u}_{\boldsymbol{g}_1}$ as in (\ref{eq2.1})-(\ref{eq2.1b}), ${\bf u}_0\in{\bf H}^{1/2}_{00}(\Gamma_0)$,
$\mbox{\bf\em g}_1\in{\bf H}^{1/2}_{00}(\Gamma_1)$, ${\rm w}_0\in H^{1/2}_{00}(\Gamma_0)$,
$ g_2\in H^{1/2}_{00}(\Gamma_3)$ and $\eta\in C^0(\mathbb{R})$ defined in (\ref{eq4}). A weak solution of (\ref{eq1.1})-(\ref{eq2.1b})
is a pair $[{\bf u},{\rm w}]\in{\bf H}_\sigma\times H^1(\Omega)$ and  $\rho=\eta(N{\bf u})$ satisfying
\begin{equation}\label{eq8}
\left\{
\begin{array}{rcl}
2\mu_1(D({\bf u}), D({\bf v}))+2\alpha\mu_1\displaystyle\int_{\Gamma_2}{\bf u}\cdot{\bf v}\,d\Gamma+\langle B[{\bf u},{\bf u},{\bf u}], {\bf v}\rangle
&\!=\!&2\mu_r({\rm rot}\, {\rm w},{\bf v})+(F({\bf u}), {\bf v})\ \forall{\bf v}\in \tilde{\bf H}_\sigma,\vspace{0.1cm}\\
\mu_2(\nabla {\rm w},\nabla {\rm z})+\langle \tilde{B}[{\bf u},{\bf u}, {\rm w}], {\rm z}\rangle+ 4\mu_r({\rm w}, {\rm z})
&\!=\!&2\mu_r({\rm rot}\,{\bf u}, {\rm z})+( G({\bf u}), {\rm z})\ \forall {\rm z}\in H^1_0(\Omega),\vspace{0.1cm}\\
\eta(N{\bf u})&\!=\!&\rho_0\quad  \mbox{ on }\Gamma_0,\\
{\bf u} &\!=\!& {\bf u}_{\boldsymbol{g}_1}\ \mbox{ on }\Gamma\setminus\Gamma_2,\\
{\rm w} &\!=\!& {\rm w}_{ g_2}\ \mbox{ on }\Gamma.
\end{array}\right.
\end{equation}
\end{definition}

In order to prove the existence of a solution to problem (\ref{eq8}), we reduce the problem to an auxiliary
problem with homogeneous boundary conditions for ${\bf u}$ on $\Gamma\setminus\Gamma_2$ and for ${\rm w}$ on $\Gamma$.
For this purpose, we introduce the following result.
\begin{lemma}\label{hopf}
Let ${\bf u}_{\boldsymbol{g}_1}$ as in (\ref{eq2.1})-(\ref{eq2.1b}), and assume that ${\bf u}_0\in{\bf H}^{1/2}_{00}(\Gamma_0)$, $\boldsymbol{g}_1\in{\bf H}^{1/2}_{00}(\Gamma_1)$. Then, for any
$\varepsilon>0,$ there exists ${\bf u}^\varepsilon\in{\bf H}_\sigma$ with ${\bf u}^\varepsilon=\boldsymbol{g}_1$
on $\Gamma_1$, ${\bf u}^\varepsilon={\bf u}_0$ on $\Gamma_0$  and ${\bf u}^\varepsilon={\bf 0}$ on $\Gamma_2$ such that
\begin{equation}\label{deshopf1}
|({\bf v} \cdot\nabla{\bf u}^\varepsilon, {\bf v})|\leq\varepsilon\|{\bf v}\|^2_{\tilde{\bf H}_\sigma}
\quad\forall {\bf v}\in\tilde{\bf H}_\sigma,\ \mbox{and}\ \Vert {\bf u}^\varepsilon\Vert_{{\bf H}_\sigma}\leq C\Vert {\bf u}_{\boldsymbol{g}_1}\Vert_{{\bf H}^{1/2}(\Gamma\setminus \Gamma_2)},
\end{equation}
where the constant $C$ depends only on $\Omega$. Moreover, if ${\rm w}_0\in H^{1/2}_{00}(\Gamma_0)$ and
$g_2\in H^{1/2}_{00}(\Gamma_3)$, then there exists ${\rm w}^e\in H^1(\Omega)$ such that
${\rm w}^e={\rm w}_0$ on $\Gamma_0$, ${\rm w}^e=g_2$ on $\Gamma_3$, and the following estimate
holds
\begin{equation}\label{deshopf2}
\|{\rm w}^e\|_{H^1}\leq C\|{\rm w}_{g_2}\|_{H^{1/2}(\Gamma)},
\end{equation}
where the constant $C$ depends only on $\Omega$.
\end{lemma}
\textit{Proof.} If ${\bf u}_0\in{\bf H}^{1/2}_{00}(\Gamma_0)$, $\boldsymbol{g}_1\in{\bf H}^{1/2}_{00}(\Gamma_1)$,
there exist $\hat{\boldsymbol g}_1\in{\bf H}^{1/2}(\Gamma),$ $\hat{\bf u}_0\in{\bf H}^{1/2}(\Gamma)$
such that $\hat{\bf u}_{0|_{\Gamma_0}}={\bf u}_{0}$, $\hat{\bf u}_{0|_{\Gamma\setminus\Gamma_0}}={\bf 0}$ ,
$\hat{\boldsymbol g}_{1|_{\Gamma\setminus\Gamma_1}}={\bf 0}$ and $\hat{\boldsymbol g}_{1|_{\Gamma_1}}=\boldsymbol{g}_1$. Thus,
$\hat{\bf u}_0+\hat{\boldsymbol g}_1\in{\bf H}^{1/2}(\Gamma)$; moreover, the integral in (\ref{eq2.1b}) implies
$\int_\Gamma(\hat{\bf u}_0+\hat{\boldsymbol g}_1)\cdot{\bf n}d\Gamma=0.$ Then, by Lemma IX.4.2 of \cite{galdi}, p. 610,
there exists ${\bf u}^\varepsilon \in {\bf H}_\sigma$ with
${\bf u}^\varepsilon=\hat{\bf u}_0+\hat{\boldsymbol g}_1$ on $\Gamma$ verifying (\ref{deshopf1}); in particular
we have that ${\bf u}^\varepsilon_{|_{\Gamma_1}}=\boldsymbol{g}_1$, ${\bf u}^\varepsilon_{|_{\Gamma_0}}={\bf u}_0$
and ${\bf u}^\varepsilon_{|_{\Gamma_2}}={\bf 0}$.
The existence of ${\rm w}^e\in H^1(\Omega)$ is well-known from the lifting theorem.
\hfill$\diamond$\vskip12pt

Rewriting $[{\bf u},{\rm w}]\in{\bf H}_\sigma\times{ H}^1(\Omega)$ in the form
${\bf u}={\bf u}^\varepsilon+\hat{\bf u}$ and ${\rm w}={\rm w}^{\rm e}+\hat{\rm w}$
with $[\hat{\bf u},\hat{\rm w}]\in\tilde{\bf H}_\sigma\times{ H}^1_0(\Omega)$
new unknown functions, from (\ref{eq8}) we obtain the following nonlinear system: Find
$[\hat{\bf u},\hat{\rm w}]\in\tilde{\bf H}_\sigma\times H^1_0(\Omega)$ such that
\begin{equation}\label{nohomo_1}
\left\{\begin{array}{l}
2\mu_1(D(\hat{\bf u}),D({\bf v}))+2\alpha\mu_1\displaystyle\int_{\Gamma_2}\hat{\bf u}\cdot{\bf v}\, d\Gamma+\langle B[\hat{\bf u}+{\bf u}^\varepsilon,\hat{\bf u}+{\bf u}^\varepsilon,\hat{\bf u}],{\bf v}\rangle
+\langle B[\hat{\bf u}+{\bf u}^\varepsilon,\hat{\bf u},{\bf u}^\varepsilon],{\bf v}\rangle\vspace{0.1cm}\\
\hspace{0.5cm}=2\mu_r({\rm rot}(\hat{\rm w}+{\rm w}^{\rm e}),{\bf v})+ (F(\hat{\bf u}+{\bf u}^\varepsilon),{\bf v})
-2\mu_1(D({\bf u}^\varepsilon),D({\bf v}))-\langle B[\hat{\bf u}+{\bf u}^\varepsilon,{\bf u}^\varepsilon,{\bf u}^\varepsilon],{\bf v}\rangle,
\vspace{0.3cm}\\
\mu_2(\nabla\hat{\rm w},\nabla{\rm z})+\langle \tilde{B}[\hat{\bf u}+{\bf u}^\varepsilon,\hat{\bf u}+{\bf u}^\varepsilon,\hat{\rm w}], {\rm z}\rangle
+\langle \tilde{B}[\hat{\bf u}+{\bf u}^\varepsilon,\hat{\bf u},{\rm w}^{\rm e}], {\rm z}\rangle+4\mu_r(\hat{\rm w},{\rm z})\vspace{0.1cm}\\
\hspace{0.5cm}=2\mu_r({\rm rot}(\hat{\bf u}+{\bf u}^\varepsilon),{\rm z})+( G(\hat{\bf u}+{\bf u}^\varepsilon),{\rm z})
-\mu_2(\nabla{\rm w}^{\rm e},\nabla{\rm z})-\langle\tilde{B}[\hat{\bf u}+{\bf u}^\varepsilon,{\bf u}^\varepsilon,{\rm w}^{\rm e}],{\rm z}\rangle
-4\mu_r({\rm w}^{\rm e},{\rm z}),
\end{array}\right.
 \end{equation}
for all $[{\bf v},{\rm z}]\in\tilde{\bf H}_\sigma\times{H}^1_0(\Omega)$.

\subsection{Statement of Boundary Control Problem}
In order to establish the statement of the boundary control problem, we suppose that
$\mathcal{U}_1\subset{\bf H}_{00}^{1/2}(\Gamma_1)$
and $\mathcal{U}_2\subset H^{1/2}_{00}(\Gamma_3)$ are nonempty sets. We consider
that ${\bf f}\in{\bf L}^2(\Omega), {\rm g}
\in {L}^2(\Omega), {\bf u}_0\in {\bf H}_{00}^{1/2}(\Gamma_0),$ ${\rm w}_0\in H^{1/2}_{00}(\Gamma_0),$
and the controls $\mbox{\bf\em g}_1\in \mathcal{U}_1$, $\mbox{\em g}_2\in \mathcal{U}_2.$ For simplicity, we denote $\mathbb{X}={\bf H}_\sigma\times{ H}^1(\Omega)\times \mathcal{U}_1\times \mathcal{U}_2$ and consider the following objective functional $J:\mathbb{X}\rightarrow \mathbb{R}$ defined by:
\begin{eqnarray}\label{objective}
J[{{\bf u}}, {\rm w}, {\mbox{\bf\em g}_1}, \mbox{\em g}_2]&=& \displaystyle\frac{\beta_1}{2}\|{\rm rot}\,{\bf u}\|^2
+\frac{\beta_2}{2}\|{\bf u}-{\bf u}_d\|^2+\frac{\beta_3}{2}\|{\rm w}-{\rm w}_d\|^2+\frac{\beta_4}{2}\|\eta(N{\bf u})-\rho_d\|^2\nonumber\\
&&+\displaystyle\frac{\beta_5}{2}\|\mbox{\bf\em g}_1\|^2_{{\bf H}^{1/2}(\Gamma_1)}
+\frac{\beta_6}{2}\|\mbox{\em g}_2\|^2_{{ H}^{1/2}(\Gamma_3)},
\end{eqnarray}
where the constants $\beta_k$ measure the cost of the control and satisfy the following conditions:
\begin{equation}\label{condBeta}
\left\{
\begin{array}{l}
i) \, \beta_i\geq 0 \mbox{ for } i=1,...,6, \mbox{ (not all zero), } \mathcal{U}_1 \subset
{\bf H}^{1/2}_{00}(\Gamma_1) \mbox{ and } \mathcal{U}_2  \subset {H}^{1/2}_{00}(\Gamma_3)
\mbox{ are bounded closed }\\
\ \ \mbox{ and convex sets};\\
ii) \, \beta_i\geq0 \mbox{ for } i=1,...,4, \, \beta_i>0 \mbox{ for } i=5,6, \,
\mathcal{U}_1  \subset {\bf H}^{1/2}_{00}(\Gamma_1) \mbox{ and } \mathcal{U}_2
\subset {H}^{1/2}_{00}(\Gamma_3)\mbox{ are closed }\\
\ \ \ \mbox{ and convex sets}.
\end{array}\right.
\end{equation}
%
In the functional (\ref{objective}), the prescribed functions ${\bf u}_d\in {\bf L}^2(\Omega),$ ${\rm w}_d\in L^2(\Omega)$ and $\rho_d\in L^2(\Omega),$ correspond to the desired states for the velocity, the microrotation velocity and the density, respectively. Then we study the following constrained minimization problem related to system (\ref{eq1.1})-(\ref{eq2.1b}):
\begin{equation}
\left\{
\begin{array}[c]{ll}
\mbox{ Find}\ [{{\bf u}}, {\rm w}, {\boldsymbol{g}_1}, \mbox{\em g}_2]\in
\mathbb{X}\ \mbox{such that, the functional}\ J[{{\bf u}}, {\rm w}, {\boldsymbol{g}_1}, \mbox{\em g}_2]\ \mbox{reaches its minimum} \\
\ \mbox{over the weak solutions of system (\ref{eq1.1})-(\ref{eq2.1b}).}
\end{array}
\right.  \label{eq1c}
\end{equation}
The set of admissible solutions of problem (\ref{eq1c}) is defined by
$$\mathcal{S}_{ad}=\{{\bf s}=[{{\bf u}}, {\rm w}, {\boldsymbol{g}_1}, \mbox{\em g}_2]\in \mathbb{X}\ \mbox{such that}\ J({\bf s})<\infty\ \mbox{and}\ {\bf s}\ \mbox{satisfies}\ (\ref{eq8})\}.$$

\section{Existence of Weak Solutions}\label{sec3}
\subsection{Linearized Problem}

For $\tilde{\bf u}\in\tilde{\bf H}_\sigma$ fixed, we consider the following linear problem: Find $\hat{\rm w}\in H^1_0(\Omega)$ such that
\begin{equation}\label{problemalineal}
\left\{\begin{array}{l}
\mu_2(\nabla\hat{\rm w},\nabla{\rm z})+\langle \tilde{B}[\tilde{\bf u}+{\bf u}^\varepsilon,\tilde{\bf u}
+{\bf u}^\varepsilon,\hat{\rm w}], {\rm z} \rangle+4\mu_r(\hat{\rm w},{\rm z})\\
\hspace{1cm}=2\mu_r({\rm rot}(\tilde{\bf u}+{\bf u}^\varepsilon),{\rm z})+(G(\tilde{\bf u}+{\bf u}^\varepsilon),{\rm z})
-\mu_2(\nabla{\rm w}^{\rm e},\nabla{\rm z})-4\mu_r({\rm w}^{\rm e}, {\rm z})\\
\hspace{1.3cm}-\langle\tilde{B}[\tilde{\bf u}+{\bf u}^\varepsilon,{\bf u}^\varepsilon,{\rm w}^{\rm e}],{\rm z}\rangle
-\langle \tilde{B}[\tilde{\bf u}+{\bf u}^\varepsilon,\tilde{\bf u},{\rm w}^{\rm e}], {\rm z}\rangle
\ \forall {\rm z}\in{ H}^1_0(\Omega),
\end{array}\right.
\end{equation}
where ${\bf u}^\varepsilon$ and ${\rm w}^e$ are given by Lemma \ref{hopf}.
For  problem (\ref{problemalineal}) we have the following result.
\begin{lemma}\label{lema3a}
If ${\rm g}\in L^2(\Omega),$ then the problem (\ref{problemalineal}) has a unique solution
$\hat{\rm w}\in H^1_0(\Omega)$. Moreover, the following inequality holds:
\begin{eqnarray}\label{desigualdad_problema_lineal}
\mu_2\|\hat{\rm w}\|_{H^1_0}&\leq& (\tilde{C}_\eta C\|{\rm w}^e\|_{H^1}
+\mu_rC)\|\tilde{\bf u}\|_{\tilde{\bf H}_\sigma}
+\tilde{C}_\eta C\|{\bf u}^\varepsilon\|_{{\bf H}_\sigma}\|{\rm w}^e\|_{H^1}
+\tilde{C}_\eta\|{\rm g}\|\nonumber\\
&&+\mu_rC(\|{\bf u}^\varepsilon\|_{{\bf H}_\sigma}
+\|{\rm w}^e\|_{H^1})+\mu_2 C\|{\rm w}^e\|_{H^1},
\end{eqnarray}
where $\tilde{C}_\eta>0$ is a constant satisfying $\|\eta(N(\tilde{\bf u}+{\bf u}^\varepsilon))\|_{\infty}
\leq\tilde{C}_\eta$, $[{\bf u}^\varepsilon, {\rm w}^e]$ is given by Lemma \ref{hopf}, and
$C>0$  is a constant that depends only on $\Omega$.
\end{lemma}
\textit{Proof.} We define the bilinear form $a_{\tilde{\bf u}}: H^1_0(\Omega)\times H^1_0(\Omega)\rightarrow\mathbb{R}$
and the linear functional ${\rm f}_{\tilde{\bf u}}: H^1_0(\Omega)\rightarrow\mathbb{R}$ by
\begin{eqnarray}
a_{\tilde{\bf u}}([\hat{\rm w},{\rm z}])&=&\mu_2(\nabla\hat{\rm w},\nabla{\rm z})+
\langle\tilde{B}[\tilde{\bf u}+{\bf u}^\varepsilon,\tilde{\bf u}+{\bf u}^\varepsilon,\hat{\rm w}],{\rm z}\rangle+4\mu_r(\hat{\rm w},{\rm z})
\ \forall {\rm z}\in{H}^1_0(\Omega),\label{eq21a}\\
\langle{\rm f}_{\tilde{\bf u}},{\rm z}\rangle&=&2\mu_r({\rm rot}(\tilde{\bf u}+{\bf u}^\varepsilon),{\rm z})+(G(\tilde{\bf u}
+{\bf u}^\varepsilon),{\rm z})-\mu_2(\nabla{\rm w}^{\rm e},\nabla{\rm z})-4\mu_r({\rm w}^{\rm e},{\rm z})\nonumber\\
&&-\langle\tilde{B}[\tilde{\bf u}+{\bf u}^\varepsilon,{\bf u}^\varepsilon,{\rm w}^{\rm e}],{\rm z}\rangle
-\langle \tilde{B}[\tilde{\bf u}+{\bf u}^\varepsilon,\tilde{\bf u},{\rm w}^{\rm e}],{\rm z}\rangle
\ \forall {\rm z}\in H^1_0(\Omega).\label{eq22a}
\end{eqnarray}
Then, from (\ref{eq21a})-(\ref{eq22a}), problem (\ref{problemalineal}) is equivalent to find $\hat{\rm w}\in H^1_0(\Omega)$ such that
\begin{equation}\label{eq11}
a_{\tilde{\bf u}}([\hat{\rm w},{\rm z}])=\langle{\rm f}_{\tilde{\bf u}},{\rm z}\rangle \ \forall {\rm z}\in H^1_0(\Omega).
\end{equation}
The bilinear form $a_{\tilde{\bf u}}$ is continuous and coercive on ${ H}^1_0(\Omega)\times{ H}^1_0(\Omega)$,
and the linear functional ${\rm f}_{\tilde{\bf u}}$ is continuous on $H^1_0(\Omega).$
Then, by Lax-Milgram Theorem, it follows that there is a unique $\hat{\rm w}\in H^1_0(\Omega)$ such that
(\ref{eq11}) is satisfied, and therefore the problem (\ref{problemalineal}) has a unique solution.

Now, in order to obtain inequality (\ref{desigualdad_problema_lineal}), by replacing ${\rm z}=\hat{\rm w}$ in (\ref{problemalineal})
and taking into account (\ref{eq6}), we have
\begin{eqnarray}\label{eq9}
\mu_2\|\nabla\hat{\rm w}\|^2+4\mu_r\|\hat{\rm w}\|^2&\leq&2\mu_r|({\rm rot}(\tilde{\bf u}+{\bf u}^\varepsilon),\hat{\rm w})|
+|(G(\hat{\bf u}+{\bf u}^\varepsilon),\hat{\rm w})|+\mu_2|(\nabla{\rm w}^{\rm e},\nabla\hat{\rm w})|\nonumber\\
&&+4\mu_r|({\rm w}^{\rm e},\hat{\rm w})|+|\langle\tilde{B}[\tilde{\bf u}
+{\bf u}^\varepsilon,{\bf u}^\varepsilon,{\rm w}^{\rm e}],\hat{\rm w}\rangle|\nonumber\\
&&+|\langle\tilde{B}[\tilde{\bf u}+{\bf u}^\varepsilon,\tilde{\bf u},{\rm w}^{\rm e}],
\hat{\rm w}\rangle|.
\end{eqnarray}
Now, we find estimates for the terms on the right hand side of (\ref{eq9}). By applying the H\"{o}lder and Poincar\'e inequalities  we get
\begin{eqnarray}
2\mu_r|({\rm rot}(\tilde{\bf u}+{\bf u}^\varepsilon),\hat{\rm w})|&\leq &
2\sqrt2\mu_r C\|\nabla(\tilde{\bf u}+{\bf u}^\varepsilon)\|\|\hat{\rm w}\|_{H^1_0}
\leq\mu_rC(\|\tilde{\bf u}\|_{\tilde{\bf H}_\sigma}+\|{\bf u}^\varepsilon\|_{{\bf H}_\sigma})\|\hat{\rm w}\|_{H^1_0},\label{eq10}\\
|(G(\tilde{\bf u}+{\bf u}^\varepsilon),\hat{\rm w})|&\leq&
\|\eta(N(\tilde{\bf u}+{\bf u}^\varepsilon))\|_{\infty}\|{\rm g}\|\|\hat{\rm w}\|_{H^1_0}
\leq\tilde{C}_\eta\|{\rm g}\|\|\hat{\rm w}\|_{H^1_0},\label{eq12}\\
\mu_2|(\nabla {\rm w}^e,\nabla\hat{\rm w})|&\leq&\mu_2\|\nabla {\rm w}^e\|\|\nabla\hat{\rm w}\|
\leq\mu_2\|{\rm w}^e\|_{H^1}\|\hat{\rm w}\|_{H^1_0},\\
4\mu_r|({\rm w}^e,\hat{\rm w})|&\leq&4\mu_r\|{\rm w}^e\|\|\hat{\rm w}\|\leq 4\mu_r C\|{\rm w}^e\|_{H^1}\|\nabla\hat{\rm w}\|
\leq\mu_r C\|{\rm w}^e\|_{H^1}\|\hat{\rm w}\|_{H^1_0},\label{eq13}\\
|\langle\tilde{B}[\tilde{\bf u}+{\bf u}^\varepsilon,{\bf u}^\varepsilon, {\rm w}^e],\hat{\rm w}\rangle|
&\leq &\|\eta(N(\tilde{\bf u}+{\bf u}^\varepsilon))\|_{\infty}\|{\bf u}^\varepsilon\|_3\|\nabla {\rm w}^e\|\|\hat{\rm w}\|_6
\leq \tilde{C}_\eta C\|\nabla{\bf u}^\varepsilon\|\|\nabla {\rm w}^e\|\|\nabla\hat{\rm w}\|\nonumber\\
&\leq& \tilde{C}_\eta C\|{\bf u}^\varepsilon\|_{{\bf H}_\sigma}\|{\rm w}^e\|_{H^1}\|\hat{\rm w}\|_{H^1_0},\label{eq14}\\
|\langle\tilde{B}[\tilde{\bf u}+{\bf u}^\varepsilon, \tilde{\bf u}, {\rm w}^e],\hat{\rm w}\rangle|
&\leq &\|\eta(N(\tilde{\bf u}+{\bf u}^\varepsilon))\|_{\infty}\|\tilde{\bf u}\|_3\|\nabla {\rm w}^e\|\|\hat{\rm w}\|_6
\leq \tilde{C}_\eta C\|\nabla\tilde{\bf u}\|\|\nabla {\rm w}^e\|\|\nabla\hat{\rm w}\|\nonumber\\
&\leq& \tilde{C}_\eta C\|\tilde{\bf u}\|_{\tilde{\bf H}_\sigma}\|{\rm w}^e\|_{H^1}\|\hat{\rm w}\|_{H^1_0}.
\label{eq15}
\end{eqnarray}
By substituting (\ref{eq10})-(\ref{eq15}) in (\ref{eq9}) we obtain
\begin{eqnarray*}
\mu_2\|\hat{\rm w}\|^2_{H^1_0}&\leq &\tilde{C}_\eta C(\|\tilde{\bf u}\|_{\tilde{\bf H}_\sigma}
+\|{\bf u}^\varepsilon\|_{{\bf H}_\sigma})\|{\rm w}^e\|_{H^1}\|\hat{\rm w}\|_{H^1_0}
+\tilde{C}_\eta\|{\rm g}\|\|\hat{\rm w}\|_{H^1_0}\\
&&+\mu_r C\big(\|\tilde{\bf u}\|_{\tilde{\bf H}_\sigma}
+\|{\bf u}^\varepsilon\|_{{\bf H}_\sigma}+\|{\rm w}^e\|_{H^1}\big)\|\hat{\rm w}\|_{H^1_0}
+\mu_2C \|{\rm w}^e\|_{H^1}\|\hat{\rm w}\|_{H^1_0},
\end{eqnarray*}
that is,
\begin{eqnarray}
\mu_2\|\hat{\rm w}\|_{H^1_0}&\leq& \tilde{C}_\eta C(\|\tilde{\bf u}\|_{\tilde{\bf H}_\sigma}
+\|{\bf u}^\varepsilon\|_{{\bf H}_\sigma})\|{\rm w}^e\|_{H^1}
+\tilde{C}_\eta\|{\rm g}\|\nonumber\\
&&+\mu_r C\big(\|\tilde{\bf u}\|_{\tilde{\bf H}_\sigma}
+\|{\bf u}^\varepsilon\|_{{\bf H}_\sigma}+\|{\rm w}^e\|_{H^1}\big)+\mu_2C \|{\rm w}^e\|_{H^1},\label{eq16}
\end{eqnarray}
which implies (\ref{desigualdad_problema_lineal}). \hfill$\diamond$\vskip12pt

\subsection{Weak Solutions}

In order to prove the existence of solutions to the problem (\ref{nohomo_1}), we define the linear operator
$T:\tilde{\bf H}_\sigma\rightarrow\tilde{\bf H}_\sigma$ as follows:  For each
 $\tilde{\bf u} \in \tilde{\bf H}_\sigma,$ let
$T\tilde{\bf u}=\hat{\bf u}$, where $\hat{\bf u}$ satisfies the following system
\begin{equation}
\begin{array}{l}\label{existence1}
2\mu_1(D(\hat{\bf u}),D({\bf v}))+2\alpha\mu_1\displaystyle\int_{\Gamma_2}\hat{\bf u}\cdot{\bf v}\,d\Gamma+\langle B[\tilde{\bf u}+{\bf u}^\varepsilon,\tilde{\bf u}+{\bf u}^\varepsilon,\hat{\bf u}],{\bf v}\rangle
+\langle B[\tilde{\bf u}+{\bf u}^\varepsilon,\hat{\bf u},{\bf u}^\varepsilon],{\bf v}\rangle\\
\\
\hspace{1cm}=2\mu_r({\rm rot}(\hat{\rm w}+{\rm w}^{\rm e}),{\bf v})+(F(\tilde{\bf u}+{\bf u}^\varepsilon),{\bf v})
-2\mu_1(D({\bf u}^\varepsilon),D({\bf v}))-\langle B[\tilde{\bf u}+{\bf u}^\varepsilon,{\bf u}^\varepsilon,{\bf u}^\varepsilon],{\bf v}\rangle,
\end{array}
\end{equation}
for all ${\bf v}\in\tilde{\bf H}_\sigma,$ and $\hat{\rm w}$ being the unique solution of the linear problem
(\ref{problemalineal}).

\begin{lemma}\label{lema3}
The operator $T:\tilde{\bf H}_\sigma\rightarrow\tilde{\bf H}_\sigma$  defined by (\ref{existence1}) is compact.
\end{lemma}
\textit{Proof.} Let $\{\tilde{\bf u}^m\}_{m\geq1}\subset\tilde{\bf H}_\sigma$
a sequence weakly convergent to $\tilde{\bf u}\in\tilde{\bf H}_\sigma$. Since the
embedding $\tilde{\bf H}_\sigma\hookrightarrow{\bf L}^p(\Omega)$ is compact for $1\leq p <\infty$ we have
$$
\tilde{\bf u}^m\rightarrow\tilde{\bf u}\ \mbox{ strongly in }{\bf L}^p(\Omega).
$$
From Remark \ref{remark2} and since the embedding $H^2(\Omega)\hookrightarrow C^0(\overline{\Omega})$ is compact,
we get that, for some subsequence of $\{\tilde{\bf u}^m\}_{m\geq1}$, still denoted by $\{\tilde{\bf u}^m\}_{m\geq1}$,
it holds that $N\tilde{\bf u}^m\rightarrow N\tilde{\bf u}$ in $C^0(\overline{\Omega});$
moreover, taking into account that $\eta\in C^0(\mathbb{R})$, we have
$\eta(N\tilde{\bf u}^m)\rightarrow \eta(N\tilde{\bf u})$ strongly in $C^0(\overline{\Omega})$,
which implies that there exists a constant $C_\eta>0$, independent of $m$, such that
\begin{equation}\label{eq_eta}
\|\eta(N(\tilde{\bf u}^m+{\bf u}^\varepsilon))\|_{\infty}\leq C_\eta.
\end{equation}
Then, denoting
$\hat{\bf u}^m=T(\tilde{\bf u}^m)$, for all  ${\bf v}\in\tilde{\bf H}_\sigma$ we get
\begin{equation}\label{existence3}
\begin{array}{l}
2\mu_1(D(\hat{\bf u}^m),D({\bf v}))+2\alpha\mu_1\displaystyle\int_{\Gamma_2}\hat{\bf u}^m\cdot{\bf v}\,d\Gamma+\langle B[\tilde{\bf u}^m+{\bf u}^\varepsilon,\tilde{\bf u}^m+{\bf u}^\varepsilon,\hat{\bf u}^m],{\bf v}\rangle
+\langle B[\tilde{\bf u}^m+{\bf u}^\varepsilon,\hat{\bf u}^m,{\bf u}^\varepsilon],{\bf v}\rangle\\
\hspace{1cm}=2\mu_r({\rm rot}(\hat{\rm w}+{\rm w}^{\rm e}),{\bf v})+(F(\tilde{\bf u}^m+{\bf u}^\varepsilon),{\bf v})
-2\mu_1(D({\bf u}^\varepsilon),D({\bf v}))-\langle B[\tilde{\bf u}^m+{\bf u}^\varepsilon,{\bf u}^\varepsilon,{\bf u}^\varepsilon],{\bf v}\rangle.
\end{array}
\end{equation}
Taking the difference between (\ref{existence3}) and (\ref{existence1}), we have
\begin{eqnarray}
&&2\mu_1(D(\hat{\bf u}^m-\hat{\bf u}),D({\bf v}))+2\alpha\mu_1\displaystyle\int_{\Gamma_2}(\hat{\bf u}^m-\hat{\bf u})\cdot{\bf v}\,d\Gamma
+\langle B[\tilde{\bf u}^m+{\bf u}^\varepsilon,\tilde{\bf u}^m
+{\bf u}^\varepsilon,\hat{\bf u}^m-\hat{\bf u}],{\bf v} \rangle\nonumber\\
&&\hspace{0.5cm}+\langle B[\tilde{\bf u}^m+{\bf u}^\varepsilon,\hat{\bf u}^m-\hat{\bf u},{\bf u}^\varepsilon],{\bf v} \rangle\nonumber\\
&&=\langle B[\tilde{\bf u}+{\bf u}^\varepsilon,\tilde{\bf u}+{\bf u}^\varepsilon,\hat{\bf u}],{\bf v} \rangle
-\langle B[\tilde{\bf u}^m+{\bf u}^\varepsilon,\tilde{\bf u}^m+{\bf u}^\varepsilon,\hat{\bf u}],{\bf v} \rangle
+\langle B[\tilde{\bf u}+{\bf u}^\varepsilon,\hat{\bf u}+{\bf u}^\varepsilon,{\bf u}^\varepsilon], {\bf v} \rangle\nonumber\\
&&\hspace{0.5cm}-\langle B[\tilde{\bf u}^m+{\bf u}^\varepsilon,\hat{\bf u}+{\bf u}^\varepsilon,{\bf u}^\varepsilon],{\bf v} \rangle
+(F(\tilde{\bf u}^m+{\bf u}^\varepsilon)-F(\tilde{\bf u}+{\bf u}^\varepsilon),{\bf v}).\label{aux1}
\end{eqnarray}
Replacing ${\bf v}=\hat{\bf u}^m-\hat{\bf u}$ in (\ref{aux1}) and using (\ref{eq6}), we deduce
\begin{eqnarray}\label{eq17}
2\mu_1\|\hat{\bf u}^m-\hat{\bf u}\|^2_{\tilde{\bf H}_\sigma}&\leq&
|\langle B[\tilde{\bf u}^m+{\bf u}^\varepsilon,\hat{\bf u}^m-\hat{\bf u},{\bf u}^\varepsilon],\hat{\bf u}^m
-\hat{\bf u}\rangle|\nonumber\\
&&+|\langle B[\tilde{\bf u}+{\bf u}^\varepsilon,\tilde{\bf u}+{\bf u}^\varepsilon,\hat{\bf u}]
-B[\tilde{\bf u}^m+{\bf u}^\varepsilon,\tilde{\bf u}^m+{\bf u}^\varepsilon,\hat{\bf u}], \hat{\bf u}^m-\hat{\bf u}\rangle|\nonumber\\
&&+|\langle B[\tilde{\bf u}+{\bf u}^\varepsilon,\hat{\bf u}+{\bf u}^\varepsilon,{\bf u}^\varepsilon]
-B[\tilde{\bf u}^m+{\bf u}^\varepsilon,\hat{\bf u}+{\bf u}^\varepsilon,{\bf u}^\varepsilon], \hat{\bf u}^m-\hat{\bf u}\rangle|\nonumber\\
&&+|(F(\tilde{\bf u}^m+{\bf u}^\varepsilon)-F(\tilde{\bf u}+{\bf u}^\varepsilon), \hat{\bf u}^m-\hat{\bf u})|.
\end{eqnarray}
Now we will bound  the terms on the right hand side of (\ref{eq17}). Using
(\ref{deshopf1}) and (\ref{eq_eta}) we obtain
\begin{eqnarray}\label{eq18}
|\langle B[\tilde{\bf u}^m+{\bf u}^\varepsilon,\hat{\bf u}^m-\hat{\bf u},{\bf u}^\varepsilon],\hat{\bf u}^m-\hat{\bf u}\rangle|
&\leq & \|\eta(N(\tilde{\bf u}^m+{\bf u}^\varepsilon))\|_{\infty}\,|\langle(\hat{\bf u}^m-\hat{\bf u})\cdot \nabla{\bf u}^\varepsilon,
\hat{\bf u}^m-\hat{\bf u}\rangle|\nonumber\\
&\leq & C_\eta\,\varepsilon \|\hat{\bf u}^m-\hat{\bf u}\|^2_{\tilde{\bf H}_\sigma}.
\end{eqnarray}
Now, notice that
\begin{eqnarray*}
&& B[\tilde{\bf u}+{\bf u}^\varepsilon,\tilde{\bf u}+{\bf u}^\varepsilon,\hat{\bf u}]
-B[\tilde{\bf u}^m+{\bf u}^\varepsilon,\tilde{\bf u}^m+{\bf u}^\varepsilon,\hat{\bf u}]\\
&&\hspace*{.5cm}= (\eta(N(\tilde{\bf u}+{\bf u}^\varepsilon))-\eta(N(\tilde{\bf u}^m+{\bf u}^\varepsilon)))
(\tilde{\bf u}+{\bf u}^\varepsilon)\cdot\nabla\hat{\bf u}
-\eta(N(\tilde{\bf u}^m+{\bf u}^\varepsilon))(\tilde{\bf u}^m-\tilde{\bf u})\cdot \nabla\hat{\bf u}.
\end{eqnarray*}
Then, by the H\"older inequality and (\ref{eq_eta}) we obtain
\begin{eqnarray}\label{eq19}
&&|\langle B[\tilde{\bf u}+{\bf u}^\varepsilon,\tilde{\bf u}+{\bf u}^\varepsilon,\hat{\bf u}]
-B[\tilde{\bf u}^m+{\bf u}^\varepsilon,\tilde{\bf u}^m+{\bf u}^\varepsilon,\hat{\bf u}], \hat{\bf u}^m-\hat{\bf u}\rangle|\nonumber\\
&&\hspace*{.5cm}\leq\|\eta(N(\tilde{\bf u}+{\bf u}^\varepsilon))-\eta(N(\tilde{\bf u}^m+{\bf u}^\varepsilon))\|_{\infty}
\,|\langle(\tilde{\bf u}+{\bf u}^\varepsilon)\cdot \nabla\hat{\bf u}, \hat{\bf u}^m-\hat{\bf u}\rangle|\nonumber\\
&&\hspace*{1cm}+\|\eta(N(\tilde{\bf u}^m+{\bf u}^\varepsilon))\|_{\infty}\,
|\langle(\tilde{\bf u}^m-\tilde{\bf u})\cdot \nabla\hat{\bf u}, \hat{\bf u}^m-\hat{\bf u}\rangle|\nonumber\\
&&\hspace*{.5cm}\leq\|\eta(N(\tilde{\bf u}+{\bf u}^\varepsilon))-\eta(N(\tilde{\bf u}^m+{\bf u}^\varepsilon))\|_{\infty}
\,\|(\tilde{\bf u}+{\bf u}^\varepsilon)\cdot \nabla\hat{\bf u}\|_{\tilde{\bf H}'}\|\hat{\bf u}^m-\hat{\bf u}\|_{\tilde{\bf H}_\sigma}\nonumber\\
&&\hspace*{1cm}+C_\eta\,\|\tilde{\bf u}^m-\tilde{\bf u}\|_3\|\hat{\bf u}\|_{\tilde{\bf H}_\sigma} \|\hat{\bf u}^m-\hat{\bf u}\|_6\nonumber\\
&&\hspace*{.5cm}\leq C\|\eta(N(\tilde{\bf u}+{\bf u}^\varepsilon))-\eta(N(\tilde{\bf u}^m+{\bf u}^\varepsilon))\|_{\infty}
\|\hat{\bf u}^m-\hat{\bf u}\|_{\tilde{\bf H}_\sigma}\nonumber\\
&&\hspace*{1cm}+CC_\eta\|\tilde{\bf u}^m-\tilde{\bf u}\|_3\|\hat{\bf u}^m-\hat{\bf u}\|_{\tilde{\bf H}_\sigma}.\hspace{1cm}
\end{eqnarray}
Again, by using the H\"{o}lder inequality we have
\begin{eqnarray}
&&|\langle B[\tilde{\bf u}+{\bf u}^\varepsilon,\hat{\bf u}+{\bf u}^\varepsilon,{\bf u}^\varepsilon]
-B[\tilde{\bf u}^m+{\bf u}^\varepsilon,\hat{\bf u}+{\bf u}^\varepsilon,{\bf u}^\varepsilon], \hat{\bf u}^m-\hat{\bf u}\rangle|\nonumber\\
&&\hspace*{1cm}\leq \|\eta(N(\tilde{\bf u}+{\bf u}^\varepsilon))-\eta(N(\tilde{\bf u}^m+{\bf u}^\varepsilon))\|_{\infty}
\|(\hat{\bf u}+{\bf u}^\varepsilon)\cdot \nabla {\bf u}^\varepsilon\|_{\tilde{\bf H}'}\|\hat{\bf u}^m-\hat{\bf u}\|_{\tilde{\bf H}_\sigma}\nonumber\\
&&\hspace*{1cm}\leq C\|\eta(N(\tilde{\bf u}+{\bf u}^\varepsilon))-\eta(N(\tilde{\bf u}^m+{\bf u}^\varepsilon))\|_{\infty}
\,\|\hat{\bf u}^m-\hat{\bf u}\|_{\tilde{\bf H}_\sigma},\label{eq20}
\end{eqnarray}
and
\begin{eqnarray}
&&|( F(\tilde{\bf u}^m+{\bf u}^\varepsilon)-F(\tilde{\bf u}+{\bf u}^\varepsilon), \hat{\bf u}^m-\hat{\bf u})|\nonumber\\
&&\hspace*{1cm}\leq \|\eta(N(\tilde{\bf u}^m+{\bf u}^\varepsilon))-\eta(N(\tilde{\bf u}+{\bf u}^\varepsilon))\|_{\infty}
\|{\bf f}\|\|\hat{\bf u}^m-\hat{\bf u}\|_{\tilde{\bf H}_\sigma}\nonumber\\
&&\hspace*{1cm}\leq C\|\eta(N(\tilde{\bf u}^m+{\bf u}^\varepsilon))-\eta(N(\tilde{\bf u}+{\bf u}^\varepsilon))\|_{\infty}
\|\hat{\bf u}^m-\hat{\bf u}\|_{\tilde{\bf H}_\sigma}.\label{eq21}
\end{eqnarray}
Replacing the inequalities (\ref{eq18})-(\ref{eq21}) in (\ref{eq17}), and taking into account (\ref{eq_eta}), we obtain

\begin{eqnarray*}
(2\mu_1-\varepsilon C_\eta)\|\hat{\bf u}^m-\hat{\bf u}\|^2_{\tilde{\bf H}_\sigma}&\leq&
C\|\eta(N(\tilde{\bf u}+{\bf u}^\varepsilon))-\eta(N(\tilde{\bf u}^m+{\bf u}^\varepsilon))\|_{\infty}
\|\hat{\bf u}^m-\hat{\bf u}\|_{\tilde{\bf H}_\sigma}\nonumber\\
&&+CC_\eta\|\tilde{\bf u}^m-\tilde{\bf u}\|_3\|\hat{\bf u}^m-\hat{\bf u}\|_{\tilde{\bf H}_\sigma}.
\end{eqnarray*}
Thus, for $\varepsilon$ small enough such that $2\mu_1-\varepsilon C_\eta>0$ we deduce that
\begin{equation}\label{eq22}
\|\hat{\bf u}^m-\hat{\bf u}\|_{\tilde{\bf H}_\sigma}\leq C(\|\eta(N(\tilde{\bf u}+{\bf u}^\varepsilon))-\eta(N(\tilde{\bf u}^m+{\bf u}^\varepsilon))\|_{\infty}
+C_\eta\|\tilde{\bf u}^m-\tilde{\bf u}\|_3).
\end{equation}
Passing to the limit in (\ref{eq22}), when $m\rightarrow\infty$, and  considering the strong convergences of
$\{\tilde{\bf u}^m\}_{m\geq1}$ in ${\bf L}^3(\Omega)$, and $\{\eta(N\tilde{\bf u}^m)\}_{m\geq1}$ in $C^0(\overline{\Omega})$, we have
$$
\|T\tilde{\bf u}^m-T\tilde{\bf u}\|_{\tilde{\bf H}_\sigma}=\|\hat{\bf u}^m-\hat{\bf u}\|_{\tilde{\bf H}_\sigma}\rightarrow 0,\ \mbox{as}\ m\rightarrow \infty.
$$
Thus, we conclude that $T$ is a compact operator. \hfill$\diamond$\vskip12pt

\begin{lemma} \label{lema4} Let $T$ be the operator defined by (\ref{existence1}), and consider the set

\begin{equation}\label{existence6}
M_\lambda=\{\tilde{\bf u}\in \tilde{\bf H}_\sigma : \tilde{\bf u}=\lambda T\tilde{\bf u} \ for \ some \ \lambda \in [0, 1]\}.
\end{equation}
If $\mu_1$ and $\mu_2$ are large enough such that
\begin{equation}\label{cotaradio}
\min\{2\mu_1-\tilde{C}_\eta C\|{\rm w}_{g_2}\|_{H^{1/2}(\Gamma)}-\mu_rC,\mu_2-\mu_rC\}>0,
\end{equation}
where $\tilde{C}_\eta$ is the constant defined in Lemma \ref{lema3a} and $C>0$
is a constant that depends only on $\Omega$, then, the set $M_\lambda$ is bounded in $\tilde{\bf H}_\sigma$.
Moreover, for $\lambda\in[0,1]$, all functions $\tilde{\bf u}=\tilde{\bf u}(\lambda)\in M_\lambda$ are contained,
independently of $\lambda$, in the open ball $B({\bf 0};R)\subset\tilde{\bf H}_\sigma,$ with
\begin{eqnarray}\label{eq23}
R&=&\tilde{C}_\eta C_1 (\|{\bf u}_{{\boldsymbol g}_1}\|_{{\bf H}^{1/2}(\Gamma\setminus\Gamma_2)}\|{\rm w}_{g_2}\|_{H^{1/2}(\Gamma)}
+ \|{\bf u}_{{\boldsymbol g}_1}\|^2_{{\bf H}^{1/2}(\Gamma\setminus\Gamma_2)}+\|{\bf f}\|
+\|{\rm g}\|) \nonumber\\
&&+C_1(\|{\bf u}_{{\boldsymbol g}_1}\|_{{\bf H}^{1/2}(\Gamma\setminus\Gamma_2)}+\|{\rm w}_{g_2}\|_{H^{1/2}(\Gamma)}),
\end{eqnarray}
where the constant $C_1>0$ depends only on $\Omega$, $\mu_1$ and $\mu_2$.

\end{lemma}
\textit{Proof.} Assuming $\lambda>0,$ for any $\tilde{\bf u}\in M_\lambda$ we can write $T\tilde{\bf u}=\frac1\lambda\tilde{\bf u};$
then, substituting $\hat{\bf u}=\frac1\lambda\tilde{\bf u}$ and
${\bf v}=\lambda\tilde{\bf u}$ in (\ref{existence1}), and taking into account (\ref{eq6}), we obtain
\begin{eqnarray}\label{eq24}
2\mu_1\|\tilde{\bf u}\|_{\tilde{\bf H}_\sigma}^2&\leq&|\langle B[\tilde{\bf u}+{\bf u}^\varepsilon,\tilde{\bf u},{\bf u}^\varepsilon],\tilde{\bf u}\rangle|
+2\mu_r\lambda|({\rm rot}(\hat{\rm w}+{\rm w}^{\rm e}),\tilde{\bf u})|+\lambda|(F(\tilde{\bf u}+{\bf u}^\varepsilon),\tilde{\bf u})|\nonumber\\
&&+2\mu_1\lambda|(D({\bf u}^\varepsilon),D(\tilde{\bf u}))|+\lambda|\langle B[\tilde{\bf u}+{\bf u}^\varepsilon,{\bf u}^\varepsilon,{\bf u}^\varepsilon],\tilde{\bf u}\rangle|.
\end{eqnarray}
Now we will bound the terms on the right hand side  of (\ref{eq24}). Using the H\"{o}lder inequality and
(\ref{deshopf1}) we obtain
\begin{eqnarray}\label{eq25}
|\langle B[\tilde{\bf u}+{\bf u}^\varepsilon,\tilde{\bf u},{\bf u}^\varepsilon],\tilde{\bf u}\rangle|
&\leq&\varepsilon\|\tilde{\bf u}\|^2_{\tilde{\bf H}_\sigma}
\|\eta(N(\tilde{\bf u}+{\bf u}^\varepsilon))\|_{\infty}\leq \varepsilon\tilde{C}_\eta \|\tilde{\bf u}\|^2_{\tilde{\bf H}_\sigma}.
\end{eqnarray}
Applying the H\"{o}lder and Poincar\'e inequalities we obtain

\begin{eqnarray}
2\mu_r\lambda|({\rm rot}(\hat{\rm w}+{\rm w}^e),\tilde{\bf u})|&\leq&
2\mu_r\lambda\|{\rm rot}(\hat{\rm w}+{\rm w}^e)\|\|\tilde{\bf u}\|
\leq\lambda\mu_r C\|\nabla(\hat{\rm w}+{\rm w}^e)\|\|\tilde{\bf u}\|_{\tilde{\bf H}_\sigma}\nonumber\\
&\leq&\lambda \mu_r C(\|\hat{\rm w}\|_{H^1_0}+\|{\rm w}^e\|_{H^1})\|\tilde{\bf u}\|_{\tilde{\bf H}_\sigma},\label{eq26}\\
\lambda|(F(\tilde{\bf u}+{\bf u}^\varepsilon),\tilde{\bf u})|
&\leq&\lambda\|\eta(N(\tilde{\bf u}+{\bf u}^\varepsilon))\|_{\infty}
\|{\bf f}\|\|\tilde{\bf u}\|_{\tilde{\bf H}_\sigma}
\leq \lambda\tilde{C}_\eta\|{\bf f}\|\|\tilde{\bf u}\|_{\tilde{\bf H}_\sigma},\label{eq27}\\
2\mu_1\lambda|(D({\bf u}^\varepsilon),D(\tilde{\bf u}))|&\leq&
2\mu_1\lambda\|D({\bf u}^\varepsilon)\|\|\tilde{\bf u}\|_{\tilde{\bf H}_\sigma}
\leq2\mu_1\lambda\|{\bf u}^\varepsilon\|_{{\bf H}_\sigma}\|\tilde{\bf u}\|_{\tilde{\bf H}_\sigma},\label{eq28}\\
\lambda|\langle B[\tilde{\bf u}+{\bf u}^\varepsilon,
{\bf u}^\varepsilon,{\bf u}^\varepsilon],\tilde{\bf u}\rangle|
&\leq&\lambda\|\eta(N(\tilde{\bf u}+{\bf u}^\varepsilon))\|_{\infty}
\|{\bf u}^\varepsilon\|_3\|\nabla{\bf u}^\varepsilon\|\|\tilde{\bf u}\|_6
\leq \lambda \tilde{C}_\eta C\|{\bf u}^\varepsilon\|^2_{{\bf H}_\sigma}\|\tilde{\bf u}\|_{\tilde{\bf H}_\sigma}.\label{eq29}
\end{eqnarray}
Replacing (\ref{eq25})-(\ref{eq29}) in (\ref{eq24}) and taking into account that $\lambda\leq1$, we have
\begin{eqnarray*}
(2\mu_1-\varepsilon \tilde{C}_\eta)\|\tilde{\bf u}\|_{\tilde{\bf H}_\sigma}^2&\leq&
\mu_r C(\|\hat{\rm w}\|_{H^1_0}+\|{\rm w}^e\|_{H^1})\|\tilde{\bf u}\|_{\tilde{\bf H}_\sigma}
+\tilde{C}_\eta\|{\bf f}\|\|\tilde{\bf u}\|_{\tilde{\bf H}_\sigma}
+2\mu_1\|{\bf u}^\varepsilon\|_{{\bf H}_\sigma}\|\tilde{\bf u}\|_{\tilde{\bf H}_\sigma}\nonumber\\
&&+C\tilde{C}_\eta\|{\bf u}^\varepsilon\|^2_{{\bf H}_\sigma}\|\tilde{\bf u}\|_{\tilde{\bf H}_\sigma},
\end{eqnarray*}
which implies
\begin{equation}\label{cotau}
(2\mu_1-\varepsilon \tilde{C}_\eta)\|\tilde{\bf u}\|_{\tilde{\bf H}_\sigma}\leq \mu_r C(\|\hat{\rm w}\|_{H^1_0}+\|{\rm w}^e\|_{H^1})
+\tilde{C}_\eta\|{\bf f}\|+2\mu_1\|{\bf u}^\varepsilon\|_{{\bf H}_\sigma}
+\tilde{C}_\eta C\|{\bf u}^\varepsilon\|^2_{{\bf H}_\sigma}.
\end{equation}
Adding (\ref{desigualdad_problema_lineal}) and (\ref{cotau}), and taking into account (\ref{deshopf1})-(\ref{deshopf2}),
we get
\begin{eqnarray*}
(2\mu_1-\varepsilon \tilde{C}_\eta)\|\tilde{\bf u}\|_{\tilde{\bf H}_\sigma}+\mu_2\|\hat{\rm w}\|_{H^1_0}
&\leq &(\tilde{C}_\eta C\|{\rm w}_{g_2}\|_{H^{1/2}(\Gamma)}+\mu_rC)\|\tilde{\bf u}\|_{\tilde{\bf H}_\sigma}
+\mu_r C\|\hat{\rm w}\|_{H^1_0}\\
&&+\tilde{C}_\eta C\|{\bf u}_{{\boldsymbol g}_1}\|_{{\bf H}^{1/2}(\Gamma\setminus\Gamma_2)}\|{\rm w}_{g_2}\|_{H^{1/2}(\Gamma)}\nonumber\\
&&+\mu_rC(\|{\bf u}_{{\boldsymbol g}_1}\|_{{\bf H}^{1/2}(\Gamma\setminus\Gamma_2)}+\|{\rm w}_{g_2}\|_{H^{1/2}(\Gamma)})\nonumber\\
&&+\mu_2 C\|{\rm w}_{g_2}\|_{H^{1/2}(\Gamma)}+2\mu_1C\|{\bf u}_{{\boldsymbol g}_1}\|_{{\bf H}^{1/2}(\Gamma\setminus\Gamma_2)}\nonumber\\
&&+\tilde{C}_\eta C\|{\bf u}_{{\boldsymbol g}_1}\|^2_{{\bf H}^{1/2}(\Gamma\setminus\Gamma_2)}
+\tilde{C}_\eta (\|{\bf f}\|+\|{\rm g}\|),
\end{eqnarray*}
and thus,
\begin{eqnarray}\label{eq31}
&&(2\mu_1-\varepsilon \tilde{C}_\eta-\tilde{C}_\eta C\|{\rm w}_{g_2}\|_{H^{1/2}(\Gamma)}
-\mu_r C)\|\tilde{\bf u}\|_{\tilde{\bf H}_\sigma}+(\mu_2-\mu_r C)\|\hat{\rm w}\|_{H^1_0}\nonumber\\
&&\hspace{1cm}\leq \tilde{C}_\eta C\|{\bf u}_{{\boldsymbol g}_1}\|_{{\bf H}^{1/2}(\Gamma\setminus\Gamma_2)}
\|{\rm w}_{g_2}\|_{H^{1/2}(\Gamma)}+\mu_rC(\|{\bf u}_{{\boldsymbol g}_1}\|_{{\bf H}^{1/2}(\Gamma\setminus\Gamma_2)}
+\|{\rm w}_{g_2}\|_{H^{1/2}(\Gamma)})\nonumber\\
&&\hspace{1.5cm}+\mu_2 C\|{\rm w}_{g_2}\|_{H^{1/2}(\Gamma)}
+2\mu_1C\|{\bf u}_{{\boldsymbol g}_1}\|_{{\bf H}^{1/2}(\Gamma\setminus\Gamma_2)}
+\tilde{C}_\eta C\|{\bf u}_{{\boldsymbol g}_1}\|^2_{{\bf H}^{1/2}(\Gamma\setminus\Gamma_2)}\nonumber\\
&&\hspace{1.5cm}+\tilde{C}_\eta (\|{\bf f}\|+\|{\rm g}\|).
\end{eqnarray}
By using  (\ref{cotaradio}), $\mu_2-\mu_r C >0$ and $\delta= 2\mu_1-\varepsilon \tilde{C}_\eta
-\tilde{C}_\eta C\|{\rm w}_{g_2}\|_{H^{1/2}(\Gamma)}-\mu_r C >0$ , for $\varepsilon$ small enough. Then, from
(\ref{eq31}) it follows

\begin{eqnarray}\label{eq3.3}
&&\delta\|\tilde{\bf u}\|_{\tilde{\bf H}_\sigma}\leq \tilde{C}_\eta C\|{\bf u}_{{\boldsymbol g}_1}\|_{{\bf H}^{1/2}(\Gamma\setminus\Gamma_2)}
\|{\rm w}_{g_2}\|_{H^{1/2}(\Gamma)}+\mu_rC(\|{\bf u}_{{\boldsymbol g}_1}\|_{{\bf H}^{1/2}(\Gamma\setminus\Gamma_2)}
+\|{\rm w}_{g_2}\|_{H^{1/2}(\Gamma)})\nonumber\\
&&\hspace{1.5cm}+\mu_2 C\|{\rm w}_{g_2}\|_{H^{1/2}(\Gamma)}
+2\mu_1C\|{\bf u}_{{\boldsymbol g}_1}\|_{{\bf H}^{1/2}(\Gamma\setminus\Gamma_2)}
+\tilde{C}_\eta C\|{\bf u}_{{\boldsymbol g}_1}\|^2_{{\bf H}^{1/2}(\Gamma\setminus\Gamma_2)}\nonumber\\
&&\hspace{1.5cm}+\tilde{C}_\eta (\|{\bf f}\|+\|{\rm g}\|),
\end{eqnarray}
which implies that $M_\lambda$ is bounded in $\tilde{\bf H}_\sigma$ for $\lambda>0$.
For $\lambda=0$ the result is trivial. The radius $R$ in (\ref{eq23}) follows
from (\ref{eq3.3}). \hfill$\diamond$\vskip12pt

With the previous results, we have the following theorem of existence of solutions for the system (\ref{nohomo_1}).

\begin{theorem}\label{teo_exist}
Let ${\bf f}\in{\bf L}^2(\Omega)$, ${\rm g}\in{ L}^2(\Omega)$, $\eta\in C^0(\mathbb{R})$, and $\mu_1$, $\mu_2$ satisfying
(\ref{cotaradio}). Then, the operator $T:\tilde{\bf H}_\sigma
\rightarrow\tilde{\bf H}_\sigma$ defined in (\ref{existence1}), with $\hat{\rm w}$ a solution of
(\ref{problemalineal}), has a fixed point $\hat{\bf u}\in\tilde{\bf H}_\sigma$ and the pair of
functions $[\hat{\bf u},\hat{\rm w}]$ is a solution of system (\ref{nohomo_1}). Furthermore, the
solution $[\hat{\bf u},\hat{\rm w}]$ satisfies the inequality
\begin{equation}\label{eq3.4}
\|\hat{\bf u}\|_{\tilde{\bf H}_\sigma}+\|\hat{\rm w}\|_{{ H}^1_0}\leq C \Theta,
\end{equation}
where  $\Theta=\|{\bf u}_{{\boldsymbol g}_1}\|_{{\bf H}^{1/2}(\Gamma\setminus\Gamma_2)}\|{\rm w}_{g_2}\|_{H^{1/2}(\Gamma)}
+\|{\bf u}_{{\boldsymbol g}_1}\|_{{\bf H}^{1/2}(\Gamma\setminus\Gamma_2)}+\|{\rm w}_{g_2}\|_{H^{1/2}(\Gamma)}+
\|{\bf u}_{{\boldsymbol g}_1}\|^2_{{\bf H}^{1/2}(\Gamma\setminus\Gamma_2)}+\|{\bf f}\|+\|{\rm g}\|$, and $C$ 
positive constant that depends only on $\Omega,\mu_1,\mu_2,\mu_r$,
$\tilde{C}_\eta$.
\end{theorem}
\textit{Proof.} From Lemmas \ref{lema3} and \ref{lema4}, it follows that the operator
$T$ and the set $M_\lambda$  satisfy the conditions of the Leray-Schauder theorem
(cf. Theorem 1.2.4, p. 42 of \cite{lukaszewicz}); therefore the operator $T$ has a fixed point,
that is, there exists $\hat{\bf u}\in\tilde{\bf H}_\sigma$ such that $T\hat{\bf u}=\hat{\bf u}$. Then,
by the definition of $T$ it follows that $\tilde{\bf u}= \hat{\bf u} $ which satisfies (\ref{existence1}),
and the auxiliary equation (\ref{problemalineal}). Thus, we concluded
that $ [\hat{\bf u}, \hat{\rm w}] $ is a solution of system (\ref{nohomo_1}).\newline
\\
Now, from (\ref{eq31}) with $\tilde{\bf u}=\hat{\bf u}$ we have
\begin{eqnarray}\label{eq3.5}
&&(2\mu_1-\varepsilon \tilde{C}_\eta-\tilde{C}_\eta C\|{\rm w}_{g_2}\|_{H^{1/2}(\Gamma)}-\mu_r C)
\|\hat{\bf u}\|_{\tilde{\bf H}_\sigma}+(\mu_2-\mu_r C)\|\hat{\rm w}\|_{H^1_0}\nonumber\\
&&\hspace{1cm}\leq \tilde{C}_\eta C\|{\bf u}_{{\boldsymbol g}_1}\|_{{\bf H}^{1/2}(\Gamma\setminus\Gamma_2)}
\|{\rm w}_{g_2}\|_{H^{1/2}(\Gamma)}+\mu_rC(\|{\bf u}_{{\boldsymbol g}_1}\|_{{\bf H}^{1/2}(\Gamma\setminus\Gamma_2)}
+\|{\rm w}_{g_2}\|_{H^{1/2}(\Gamma)})\nonumber\\
&&\hspace{1.5cm}+\mu_2 C\|{\rm w}_{g_2}\|_{H^{1/2}(\Gamma)}
+2\mu_1C\|{\bf u}_{{\boldsymbol g}_1}\|_{{\bf H}^{1/2}(\Gamma\setminus\Gamma_2)}
+\tilde{C}_\eta C\|{\bf u}_{{\boldsymbol g}_1}\|^2_{{\bf H}^{1/2}(\Gamma\setminus\Gamma_2)}\nonumber\\
&&\hspace{1.5cm}+\tilde{C}_\eta (\|{\bf f}\|+\|{\rm g}\|).
\end{eqnarray}
From (\ref{cotaradio}) we have that $\hat\delta=\min\{2\mu_1-\varepsilon \tilde{C}_\eta
-\tilde{C}_\eta C\|{\rm w}_{g_2}\|_{H^{1/2}(\Gamma)}-\mu_r C, \mu_2-\mu_r C\}>0$, for $\varepsilon$ small enough. Then,
from (\ref{eq3.5}) we deduce that
\begin{eqnarray*}
\hat\delta(\|\hat{\bf u}\|_{\tilde{\bf H}_\sigma}+\|\hat{\rm w}\|_{H^1_0}) &\leq &
C(\|{\bf u}_{{\boldsymbol g}_1}\|_{{\bf H}^{1/2}(\Gamma\setminus\Gamma_2)}\|{\rm w}_{g_2}\|_{H^{1/2}(\Gamma)}
+\|{\bf u}_{{\boldsymbol g}_1}\|_{{\bf H}^{1/2}(\Gamma\setminus\Gamma_2)}+\|{\rm w}_{g_2}\|_{H^{1/2}(\Gamma)}\\
&&+\|{\bf u}_{{\boldsymbol g}_1}\|^2_{{\bf H}^{1/2}(\Gamma\setminus\Gamma_2)}+\|{\bf f}\|+\|{\rm g}\|),
\end{eqnarray*}
which implies inequality (\ref{eq3.4}).\hfill$\diamond$\vskip12pt

As a consequence of Theorem \ref{teo_exist}, we have the following result.
\begin{theorem}\label{teo_exist_nohomo}
Let ${\bf f}\in{\bf L}^2(\Omega)$, ${\rm g}\in L^2(\Omega)$,
${\bf u}_0\in{\bf H}^{1/2}_{00}(\Gamma_0)$, $\boldsymbol{g}_1\in{\bf H}^{1/2}_{00}(\Gamma_1)$,
${\rm w}_0\in H^{1/2}_{00}(\Gamma_0)$, $g_2\in H^{1/2}_{00}(\Gamma_3)$, $\eta\in C^0(\mathbb{R})$, $\mu_1$ and $\mu_2$
satisfying the condition (\ref{cotaradio}) given in Theorem \ref{teo_exist}. Then, there exists functions $[{\bf u},{\rm w}]\in{\bf H}_\sigma\times H^1(\Omega)$
and $\rho=\eta(N{\bf u})$ satisfying the nonhomogeneous system (\ref{eq8}). Moreover, the solution  $[{\bf u},{\rm w}]$ satisfies the following inequality
\begin{equation}\label{eq36}
\|{\bf u}\|_{{\bf H}_\sigma}+\|{\rm w}\|_{{ H}^1}\leq C \Theta,
\end{equation}
where $C$ is a positive constant  depending on $\Omega,\mu_1,\mu_2$, and $\Theta$
is defined as in (\ref{eq3.4}).
\end{theorem}
\textit{Proof.}
Considering the solution $[\hat{\bf u}, \hat{\rm w}] \in\tilde{\bf H}_\sigma\times H^1_0(\Omega)$ of system
(\ref{nohomo_1}) given by Theorem \ref{teo_exist}, we deduce that there exists a solution
$[{\bf u},{\rm w}]\in{\bf H}_\sigma\times H^1(\Omega)$ for system (\ref{eq8}), where ${\bf u}=\hat{\bf u}+{\bf u}^\varepsilon$
and ${\rm w}=\hat{\rm w}+{\rm w}^e$. Moreover, by using triangle inequality we have
\begin{eqnarray*}
\|{\bf u}\|_{{\bf H}_\sigma}+\|{\rm w}\|_{H^1}&\leq&\|\hat{\bf u}\|_{\tilde{\bf H}_\sigma}
+\|\hat{\rm w}\|_{H^1_0}+\|{\bf u}^\varepsilon\|_{{\bf H}_\sigma}+\|{\rm w}^e\|_{H^1}\\
&\leq&C\Theta+C(\|{\bf u}_{{\boldsymbol g}_1}\|_{{\bf H}^{1/2}(\Gamma\setminus\Gamma_2)}+\|{\rm w}_{g_2}\|_{H^{1/2}(\Gamma)}),
\end{eqnarray*}
which implies inequality (\ref{eq36}).\hfill$\diamond$\vskip12pt

\section{Existence of an Optimal Solution}\label{sec4}
In this section, we will prove the existence of an optimal solution
for problem (\ref{eq1c}). We remember that the set of admissible solutions of problem (\ref{eq1c})
is defined by
$$
\mathcal{S}_{ad}=\{{\bf s}=[{\bf u},{\rm w},\boldsymbol{g}_1,g_2]\in\mathbb{X}
\mbox{ such that } J({\bf s})<\infty \mbox{ and } {\bf s} \mbox{ satisfies } (\ref{eq8})\},
$$
where $\mathbb{X}={\bf H}_\sigma\times H^1(\Omega)\times\mathcal{U}_1\times\mathcal{U}_2$. We have the following result.
\begin{theorem}\label{teorem3}
Under the conditions of Theorem \ref{teo_exist_nohomo}, if one the conditions $(i)$ or $(ii)$
in (\ref{condBeta}) is satisfied, then the problem (\ref{eq1c}) has at least one  optimal solution
$\tilde{\bf s}=[\tilde{\bf u},\tilde{\rm w},\tilde{\boldsymbol g}_1,\tilde{\mbox{\it g}}_2]\in \mathcal{S}_{ad}$.
\end{theorem}
\textit{Proof.} From Theorem \ref{teo_exist_nohomo} we have that $S_{ad}$ is
nonempty, and since the functional $J$ is bounded below, there exists  a minimizing
sequence  $\{{\bf s}^m=[{\bf u}^m, {\rm w}^m, {\boldsymbol g}_1^m,  g_2^m]\}_{m\geq 1}\in
\mathcal{S}_{ad}, m\in \mathbb{N},$ such that $\displaystyle\lim_{m\rightarrow\infty}
J({\bf s}^m)=\inf_{{\bf s}\in \mathcal{S}_{ad}}J({\bf s})$. Moreover, ${\bf s}^m$ satisfies the
system (\ref{eq8}), that is,
\begin{equation}\label{eqm}
\left\{\begin{array}{rcl}
2\mu_1(D({\bf u}^m),D({\bf v}))+2\alpha\mu_1\displaystyle\int_{\Gamma_2}{\bf u}^m\cdot{\bf v}\,d\Gamma+\langle B[{\bf u}^m,{\bf u}^m,{\bf u}^m], {\bf v}\rangle
\!&\!=\!&\!2\mu_r({\rm rot}\, {\rm w}^m,{\bf v})+(F({\bf u}^m), {\bf v}),\vspace{0.1cm}\\
\mu_2(\nabla{\rm w}^m,\nabla{\rm z})+\langle \tilde{B}[{\bf u}^m,{\bf u}^m,{\rm w}^m], {\rm z}\rangle+ 4\mu_r({\rm w}^m, {\rm z})
\!&\!=\!&\!2\mu_r({\rm rot}\,{\bf u}^m, {\rm z})+ (G({\bf u}^m), {\rm z}),\vspace{0.1cm}\\
{\bf u}^m\!&\!=\!&\!{\bf u}_{{\boldsymbol g}_1^m}\quad \mbox{ on }\Gamma\setminus\Gamma_2,\\
{\rm w}^m\!&\!=\!&\!{\rm w}_{g_2^m}\quad\mbox{ on }\Gamma,\\
\eta(N{\bf u}^m)\!&\!=\!&\!\rho_0\quad \ \ \mbox{ on }\Gamma_0,
\end{array}\right.
\end{equation}
for all $[{\bf v},{\rm z}] \in \tilde{\bf H}_\sigma\times H^1_0(\Omega),$ and $\eta\in C^0(\mathbb{R})$.\newline
If one of the conditions  $(i)$ or $(ii)$ in (\ref{condBeta}) is satisfied, we have that
there exists  a constant $C$ independent of $m$  such that
$\|{\boldsymbol g}_1^m\|^2_{{\bf H}^{1/2}(\Gamma_1)}+\|g_2^m\|^2_{H^{1/2}(\Gamma_3)}\leq C;$
consequently, (\ref{eq36}) implies $\|{\bf u}^m\|_{{\bf H}_\sigma}+\|{\rm w}^m\|_{H^1}\leq C$. Therefore, since $\mathcal{U}_1\times\mathcal{U}_2$  is closed convex subset of
${\bf H}^{1/2}_{00}(\Gamma_1)\times H^{1/2}_{00}(\Gamma_3)$,  there exists
$\tilde{\bf s}=[\tilde{\bf u},\tilde{\rm w},\tilde{\boldsymbol g}_1,\tilde{g}_2]
\in \mathbb{X}$ and some subsequence of $\{{\bf s}^m\}_{m\geq 1}$, still denoted by
$\{{\bf s}^m\}_{m\geq1}$, such that when $m\rightarrow \infty$ we have
\begin{equation}\label{formu_minim10}
\begin{array}{lll}
&&{\bf u}^m\rightarrow \tilde{\bf u}\ \mbox{weakly in}\ {\bf H}_\sigma\ \mbox{and
strongly in}\ {\bf L}^2(\Omega),\\
&&{\rm w}^m\rightarrow \tilde{\rm w}\ \mbox{weakly in}\ {H}^1(\Omega)\ \mbox{and
strongly in}\ L^2(\Omega),\\
&&\boldsymbol{g}^m_1\rightarrow \tilde{\boldsymbol{g}}_1\ \mbox{weakly in}\
{\bf H}^{1/2}_{00}(\Gamma_1)\
\mbox{and strongly in}\ {\bf L}^2(\Gamma_1),\\
&&\mbox{\em g}^m_2\rightarrow \tilde{\mbox{\em g}}_2\ \mbox{weakly in}\
{H}^{1/2}_{00}(\Gamma_3)\
\mbox{and strongly in}\ { L}^2(\Gamma_3).
\end{array}
\end{equation}
Moreover, since ${\bf u}^m={\bf u}_{\boldsymbol{g}_1^m}$ on $\Gamma\setminus\Gamma_2$,
${\rm w}^m={\rm w}_{g_2^m}$ on $\Gamma$, and $\eta(N{\bf u}^m)=\rho_0$
on $\Gamma_0$, from (\ref{formu_minim10}) it follows that $\tilde{\bf u}={\bf u}_{\tilde{\boldsymbol{g}}_1}$
on $\Gamma\setminus\Gamma_2$, $\tilde{\rm w}={\rm w}_{\tilde{g}_2}$ on $\Gamma$,
and $\eta(N\tilde{\bf u})=\rho_0$ on $\Gamma_0$;
thus, $\tilde{\bf s}$ satisfies the boundary conditions in (\ref{eq8}). A standard procedure
permits to pass to the limit in (\ref{eqm}) when $m$ goes to $\infty$, and
then, $\tilde{\bf s}$ is solution of the system (\ref{eq8}). Thus, we have
$\tilde{\bf s}\in \mathcal{S}_{ad}$ and
\begin{eqnarray}
\lim_{m\rightarrow \infty}J({\bf s}^m)=\inf_{{\bf s}\in
S_{ad}}J({\bf s})\leq J(\tilde{\bf s}).\label{eqn6}
\end{eqnarray}
Also, since the functional $J$ is weakly lower semicontinuous on  $\mathcal{S}_{ad}$,
we get that $J(\tilde{\bf s})\leq \displaystyle\lim_{m\rightarrow\infty}\inf J({\bf s}^m)$. Therefore,
from (\ref{eqn6}) and the last inequality, we conclude that $\displaystyle J[\tilde{\bf u}, \tilde{\rm w},
\tilde{\boldsymbol{g}}_1,\tilde{\mbox{\em g}}_2]=J(\tilde{\bf s})=\min_{{\bf s}\in S_{ad}}J({\bf s}),$ which
implies the existence of an optimal solution for the control problem
(\ref{eq1c}).\hfill$\diamond$\vskip12pt

\section{Necessary Optimality Conditions and an Optimality System}\label{sec5}

This section is devoted to obtain an optimality system to problem (\ref{eq1c}). We shall use
the theorem of Lagrange multipliers to turn the constrained optimization problem (\ref{eq1c})
into an unconstrained one. In order to prove the existence of Lagrange multipliers, we use a penalty method. This method consists in introducing
a family of penalized problems $(P)_{\varepsilon}$, $\varepsilon>0,$ whose solutions converge towards a solution to
problem  (\ref{eq1c}); then we derive the optimality conditions for problems $(P)_\varepsilon$ and
finally, we pass to the limit in these optimality conditions.  This method has been previously used in \cite{Arbergel, Illarionov, Lee} in the context of Navier-Stokes and Boussinesq equations.

\subsection{Penalized Problem}
For an optimal solution $\tilde{\bf s}=[\tilde{\bf u},\tilde{\rm w},\tilde{\boldsymbol g}_1,\tilde{g}_2]$
of the optimal control problem (\ref{eq1c}) we consider the following family of auxiliary
extremal problems: Find
\begin{equation}\label{nec1}
\min J_{\varepsilon}[{\bf u}, {\rm w},\boldsymbol{g}_1, {\mbox{\it g}}_2], \qquad  [{\bf u}, {\rm w},
\boldsymbol{g}_1, {\mbox{\it g}}_2]\in \mathbb{X},
\end{equation}
where for any $\varepsilon >0$ the functional
$J_\varepsilon: \mathbb{X} \rightarrow \mathbb{R}$
is defined by
\begin{eqnarray}\label{nec2}
J_{\varepsilon}[{\bf u}, {\rm w}, \boldsymbol{g}_1,{\mbox{\it g}}_2]
&=&J[{\bf u}, {\bf w},\boldsymbol{g}_1,{\mbox{\it g}}_2]
+\frac12\|{\bf u}-\tilde{\bf u}\|^2_{{\bf H}_\sigma}+\frac12\|{\rm w}-\tilde{\rm w}\|^2_{H^1}
+\frac12\|{\boldsymbol{g}}_1-\tilde{\boldsymbol{g}}_1\|^2_{{\bf H}^{1/2}(\Gamma_1)}\nonumber\\
&&+\frac12\|g_2-\tilde{\mbox{\it g}}_2\|^2_{H^{1/2}(\Gamma_3)}
+\frac{1}{2\varepsilon}\|\mu_1A{\bf u}+B[{\bf u},{\bf u},{\bf u}]-2\mu_r{\rm rot}\,{\rm w}-F({\bf u})\|^2_{\tilde{\bf H}'}\nonumber\\
&&+\frac{1}{2\varepsilon}\|\mu_2 \tilde{A} {\rm w}+\tilde{B}[{\bf u},{\bf u}, {\rm w}]+4\mu_r{\rm w}
-2\mu_r{\rm rot}\,{\bf u}-G({\bf u})\|^2_{H^{-1}}\nonumber\\
&&+\frac{1}{2\varepsilon}\|{\bf u}-{\bf u}_{{\boldsymbol{g}}_1}\|^2_{{\bf H}^{1/2}(\Gamma\setminus\Gamma_2)}+
\frac{1}{2\varepsilon}\|{\rm w}- {\rm w}_{{g}_2}\|^2_{H^{1/2}(\Gamma)}.
\end{eqnarray}
In (\ref{nec2}), $J[{\bf u}, {\rm w},{\boldsymbol{g}}_1,{\mbox{\it g}}_2]$ is the functional defined in (\ref{eq1c}), the operators
$B$, $\tilde{B}$, $F$, and $G$ are defined in (\ref{oper1}), and the operators
$A:{\bf H}_\sigma\rightarrow\tilde{\bf H}'$ and $\tilde{A}:H^1(\Omega)\rightarrow H^{-1}(\Omega)$ are defined by
\begin{equation}\label{operA}
\langle A{\bf u},{\bf v}\rangle=2(D({\bf u}),D({\bf v}))+2\alpha\int_{\Gamma_2}{\bf u}\cdot{\bf v}\, d\Gamma\quad \forall {\bf v}\in\tilde{\bf H}_\sigma,\qquad
\langle\tilde{A}{\rm w},{\rm z}\rangle=(\nabla{\rm w},\nabla{\rm z})\quad\forall{\rm z}\in H^1_0(\Omega),
\end{equation}
where $\alpha$ is given in (\ref{eq2.1}).
\begin{remark}
Since $\tilde{\bf s}$ satisfies (\ref{eq8}), then
by definition of $J_\varepsilon$ it holds
\begin{equation}\label{nec3}
J_\varepsilon(\tilde{\bf s})=J(\tilde{\bf s}).
\end{equation}
\end{remark}

Following the proof of Theorem \ref{teorem3}, and recalling that the functional
$J_\varepsilon$ is weakly lower semincontinuous, we can prove that there exists an optimal solution
of problem (\ref{nec1})-(\ref{nec2}), that is, for each $\varepsilon>0$, there exists
${\bf s}^\varepsilon=[{\bf u}^\varepsilon,{\rm w}^\varepsilon,\boldsymbol{g}_1^\varepsilon,{\mbox{\it g}}_2^\varepsilon]
\in \mathbb{X}$ such that
$$
J_\varepsilon({\bf s}^\varepsilon)=\min_{{\bf s}\in\mathbb{X}}
J_\varepsilon({\bf s}).
$$
\begin{theorem}\label{theo_conver}
Assume the hypotheses of Theorem \ref{teorem3}. For each $\varepsilon>0$, let
${\bf s}^\varepsilon=[{\bf u}^\varepsilon,{\rm w}^\varepsilon, \boldsymbol{g}_1^\varepsilon,\mbox{\it g}_2^\varepsilon]$
a solution of problem (\ref{nec1})-(\ref{nec2}). Then, as $\varepsilon\rightarrow0,$
the sequence $\{{\bf s}^\varepsilon=[{\bf u}^\varepsilon,{\rm w}^\varepsilon,\boldsymbol{g}_1^\varepsilon,\mbox{\it g}_2^\varepsilon]\}_{\varepsilon>0}$
satisfies the following convergences
\begin{eqnarray}
[{\bf u}^\varepsilon, {\rm w}^\varepsilon, \boldsymbol{g}_1^\varepsilon,{\mbox{\it g}}_2^\varepsilon]  & \longrightarrow &
[\tilde{\bf u},\tilde{\rm w},\tilde{\boldsymbol{g}}_1, \tilde{\mbox{\it g}}_2] \ \mbox{ strongly in } \ \mathcal{S}_{ad},\label{nec4}\\
J_\varepsilon[{\bf u}^\varepsilon, {\rm w}^\varepsilon,\boldsymbol{g}_1^\varepsilon, \mbox{\it g}_2^\varepsilon]
&\longrightarrow &J[\tilde{\bf u},\tilde{\rm w}, \tilde{\boldsymbol{g}}_1,\tilde{\mbox{\it g}}_2],\label{nec8}
\end{eqnarray}
where $\tilde{\bf s}=[\tilde{\bf u},\tilde{\rm w},\tilde{\boldsymbol g}_1,\tilde{\mbox{\it g}}_2]\in\mathcal{S}_{ad}$
is an optimal solution of problem (\ref{eq1c}).
\end{theorem}

\textit{Proof.} Since $\tilde{\bf s}\in\mathcal{S}_{ad}\subset\mathbb{X}$
and the functional $J_\varepsilon$ attains its minimum in $[{\bf u}^\varepsilon,{\rm w}^\varepsilon,\boldsymbol{g}_1^\varepsilon,
\mbox{\it g}_2^\varepsilon]$, we get that $J_\varepsilon({\bf s}^\varepsilon)\leq J_\varepsilon(\tilde{\bf s});$  thus,  equality
(\ref{nec3}) implies
\begin{equation}\label{nec9}
J_\varepsilon({\bf s}^\varepsilon)\leq J(\tilde{\bf s}).
\end{equation}
Observing that $J({\bf s}^\varepsilon)\geq0$, from (\ref{nec2}) we obtain
\begin{equation*}
\frac12\|{\bf u}^\varepsilon-\tilde{\bf u}\|^2_{{\bf H}_\sigma}+\frac12\|{\rm w}^\varepsilon-\tilde{\rm w}\|^2_{H^1}
+\frac12\|\boldsymbol{g}_1^\varepsilon-\tilde{\boldsymbol{g}}_1\|^2_{{\bf H}^{1/2}(\Gamma_1)}
+\frac12\|\mbox{\it g}_2^\varepsilon-\tilde{\mbox{\it g}}_2\|^2_{H^{1/2}(\Gamma_3)}
\leq J_\varepsilon({\bf s}^\varepsilon).
\end{equation*}
Therefore,
\begin{eqnarray*}
&&\|{\bf u}^\varepsilon\|_{{\bf H}_\sigma}^2+\|{\rm w}^\varepsilon\|_{H^1}^2+\|\boldsymbol{g}_1^\varepsilon\|_{{\bf H}^{1/2}(\Gamma_1)}^2
+\|\mbox{\it g}_2^\varepsilon\|_{H^{1/2}(\Gamma_3)}^2\nonumber\\
&&\hspace{1.5cm}\leq 2(\|{\bf u}^\varepsilon-\tilde{\bf u}\|^2_{{\bf H}_\sigma}+\|{\rm w}^\varepsilon-\tilde{\rm w}\|^2_{H^1}
+\|\boldsymbol{g}_1^\varepsilon-\tilde{\boldsymbol{g}}_1\|^2_{{\bf H}^{1/2}(\Gamma_1)}
+\|\mbox{\it g}_2^\varepsilon-\tilde{\mbox{\it g}}_2\|^2_{H^{1/2}(\Gamma_3)})\nonumber\\
&&\hspace{2cm} +2(\|\tilde{\bf u}\|_{{\bf H}_\sigma}^2+\|\tilde{\rm w}\|_{H^1}^2
+\|\tilde{\boldsymbol{g}}_1\|_{{\bf H}^{1/2}(\Gamma_1)}^2+\|\tilde{\mbox{\it g}}_2\|_{H^{1/2}(\Gamma_3)}^2)\nonumber\\
&&\hspace{1.5cm}\leq 4J_\varepsilon({\bf s}^\varepsilon)+2(\|\tilde{\bf u}\|_{{\bf H}_\sigma}^2+\|\tilde{\rm w}\|_{H^1}^2
+\|\tilde{\boldsymbol{g}}_1\|_{{\bf H}^{1/2}(\Gamma_1)}^2+\|\tilde{\mbox{\it g}}_2\|_{H^{1/2}(\Gamma_3)}^2),
\end{eqnarray*}
and using (\ref{nec9}) we obtain
\begin{equation}\label{nec10}
\|{\bf u}^\varepsilon\|_{{\bf H}_\sigma}^2+\|{\rm w}^\varepsilon\|_{H^1}^2+\|\boldsymbol{g}_1^\varepsilon\|_{{\bf H}^{1/2}(\Gamma_1)}^2
+\|\mbox{\it g}_2^\varepsilon\|_{H^{1/2}(\Gamma_3)}^2\leq 4 J(\tilde{\boldsymbol s})+C \leq C,
\end{equation}
where $C$ is a constant independent of $\varepsilon$.
Thus, (\ref{nec10}) implies that
there exists ${\bf s}=[{\bf u},{\rm w},\boldsymbol{g}_1,\mbox{\it g}_2]\in\mathbb{X}$
and a subsequence of $\{{\bf s}^\varepsilon=[{\bf u}^\varepsilon, {\rm w}^\varepsilon,\boldsymbol{g}_1^\varepsilon,\mbox{\it g}_2^\varepsilon ]\}_{\varepsilon>0},$ still denoted by $\{{\bf s}^\varepsilon\}_{\varepsilon>0},$
such that
\begin{equation}\label{nec11}
[{\bf u}^\varepsilon, {\rm w}^\varepsilon,\boldsymbol{g}_1^\varepsilon,\mbox{\it g}_2^\varepsilon]
\longrightarrow [{\bf u}, {\rm w},\boldsymbol{g}_1,\mbox{\it g}_2] \ \mbox{ weakly in }
\mathbb{X},\ \mbox{as}\ \varepsilon\rightarrow0.
\end{equation}
Since the embeddings
${\bf H}_\sigma\hookrightarrow{\bf L}^2(\Omega)$, $H^1(\Omega)\hookrightarrow L^2(\Omega)$,
$\mathcal{U}_1\hookrightarrow{\bf L}^2(\Gamma_1)$ and $\mathcal{U}_2\hookrightarrow L^2(\Gamma_3)$
are compact, we have
\begin{equation}\label{nec12}
[{\bf u}^\varepsilon, {\rm w}^\varepsilon,\boldsymbol{g}_1^\varepsilon,\mbox{\it g}_2^\varepsilon]
\longrightarrow [{\bf u},{\rm w},\boldsymbol{g}_1,\mbox{\it g}_2]\ \mbox{ strongly in }
{\bf L}^2(\Omega)\times L^2(\Omega)\times{\bf L}^2(\Gamma_1)\times L^2(\Gamma_3).
\end{equation}
Following the proof of Theorem \ref{teorem3}, by considering
(\ref{nec11}) and (\ref{nec12}), we obtain that ${\bf s}=[{\bf u},{\rm w},\boldsymbol{g}_1,\mbox{\it g}_2]\in\mathcal{S}_{ad}$.
Furthermore, since ${\bf s}^\varepsilon\in\mathbb{X}$,
by the definition of $J_\varepsilon$ given in (\ref{nec2}), we obtain
$$
J({\bf s}^\varepsilon)\leq J_\varepsilon({\bf s}^\varepsilon),
$$
and by using (\ref{nec9}), we get
\begin{equation}\label{nec13}
J({\bf s}^\varepsilon)\leq J_\varepsilon({\bf s}^\varepsilon)\leq J(\tilde{\bf s}).
\end{equation}
Since ${\bf s}\in\mathcal{S}_{ad}$ and
$J$ is weakly lower semicontinuous on  $\mathcal{S}_{ad}$,  from (\ref{nec11})
and (\ref{nec12}) we obtain
$$
J[{\bf u},{\rm w},\boldsymbol{g}_1,\mbox{\it g}_2]\leq\lim_{\varepsilon\rightarrow0}\inf J({\bf s}^\varepsilon)
\leq J[\tilde{\bf u},\tilde{\rm w},\tilde{\boldsymbol{g}}_1,\tilde{\mbox{\it g}}_2]
\leq J[{\bf u},{\rm w},{\boldsymbol{g}}_1,\mbox{\it g}_2],
$$
and taking into account that $\tilde{\bf s}=[\tilde{\bf u},\tilde{\rm w},\tilde{\boldsymbol{g}}_1,\tilde{\mbox{\it g}}_2]$
is an optimal solution for control problem (\ref{eq1c}), we conclude that
$[{\bf u},{\rm w},{\boldsymbol{g}}_1,\mbox{\it g}_2]=[\tilde{\bf u},\tilde{\rm w},\tilde{\boldsymbol{g}}_1,\mbox{\it g}_2]$. Then,
from (\ref{nec11})-(\ref{nec12}), as $\varepsilon\rightarrow0$, we get

\begin{eqnarray}
&&[{\bf u}^\varepsilon,{\rm w}^\varepsilon,{\boldsymbol{g}}_1^\varepsilon,\mbox{\it g}_2^\varepsilon]\longrightarrow
[\tilde{\bf u},\tilde{\rm w},\tilde{\boldsymbol{g}}_1,\tilde{\mbox{\it g}}_2]\ \mbox{ weakly in }\mathcal{S}_{ad}
\subset \mathbb{X},\label{nec14}\\
&&[{\bf u}^\varepsilon, {\rm w}^\varepsilon,\boldsymbol{g}_1^\varepsilon,\mbox{\it g}_2^\varepsilon ]
\longrightarrow [\tilde{\bf u},\tilde{\rm w},\tilde{\boldsymbol g}_1,\tilde{\mbox{\it g}}_2]\ \mbox{ strongly in }
{\bf L}^2(\Omega)\times L^2(\Omega)\times{\bf L}^2(\Gamma_1)\times L^2(\Gamma_3).\qquad \label{nec14a}
\end{eqnarray}
Now, we observe that
\begin{eqnarray}
\|{\bf u}^\varepsilon\|^2_{{\bf H}_\sigma}-\|\tilde{\bf u}\|^2_{{\bf H}_\sigma}
&=& \|{\bf u}^\varepsilon\|^2-\|\tilde{\bf u}\|^2+ (\nabla {\bf u}^\varepsilon,
\nabla {\bf u}^\varepsilon-\nabla \tilde{\bf u})
+(\nabla \tilde{\bf u}, \nabla {\bf u}^\varepsilon-\nabla \tilde{\bf u}),\label{nec15}\\
\|{\rm w}^\varepsilon\|^2_{H^1}-\|\tilde{\rm w}\|^2_{H^1}&=&\|{\rm w}^\varepsilon\|^2
-\|\tilde{\rm w}\|^2+ (\nabla {\rm w}^\varepsilon,
\nabla{\rm w}^\varepsilon-\nabla \tilde{\rm w})+(\nabla\tilde{\rm w}, \nabla{\rm w}^\varepsilon
-\nabla \tilde{\rm w}),\label{nec16}
\end{eqnarray}
where $\nabla{\bf u}^\varepsilon,\ \nabla\tilde{\bf u}\in{\bf L}^2(\Omega)$ and
$\nabla{\rm w}^\varepsilon,\,\nabla\tilde{\rm w}\in L^2(\Omega)$. Therefore, from (\ref{nec14})-(\ref{nec14a}), as $\varepsilon\rightarrow0$, we get
\begin{eqnarray}
(\nabla{\bf u}^\varepsilon,\nabla{\bf u}^\varepsilon-\nabla\tilde{\bf u})
+(\nabla\tilde{\bf u},\nabla{\bf u}^\varepsilon-\nabla\tilde{\bf u})\longrightarrow0
&\mbox{ and }&\|{\bf u}^\varepsilon\|\longrightarrow\|\tilde{\bf u}\|,\label{nec19}\\
(\nabla{\rm w}^\varepsilon,\nabla{\rm w}^\varepsilon-\nabla\tilde{\rm w})
+(\nabla\tilde{\rm w},\nabla{\rm w}^\varepsilon-\nabla\tilde{\rm w})\longrightarrow0
&\mbox{ and }&\|{\rm w}^\varepsilon\|\longrightarrow\|\tilde{\rm w}\|,\label{nec20}
\end{eqnarray}
which, together with (\ref{nec15})-(\ref{nec16}), implies that
$$
\|{\bf u}^\varepsilon\|^2_{{\bf H}_\sigma}\longrightarrow \|\tilde{\bf u}\|^2_{{\bf H}_\sigma}\ \mbox{ and }\
\|{\rm w}^\varepsilon\|^2_{{ H}^1}\longrightarrow \|\tilde{\rm w}\|^2_{{H}^1}.
$$
In particular,
\begin{equation}\label{nec21}
\limsup_{\varepsilon\rightarrow0}\|{\bf u}^\varepsilon\|_{{\bf H}_\sigma}=\|\tilde{\bf u}\|_{{\bf H}_\sigma}\ \mbox{ and }\
\limsup_{\varepsilon\rightarrow0}\|{\rm w}^\varepsilon\|_{{H}^1}=\|\tilde{\rm w}\|_{{H}^1}.
\end{equation}
Since ${\bf H}_\sigma\times H^1(\Omega)$ is a Hilbert space, it is
uniformly convex; then from (\ref{nec21}) and Proposition 3.32 of \cite{brezis}
p. 78, as $\varepsilon\rightarrow0,$ we obtain
\begin{equation}\label{strong1}
[{\bf u}^\varepsilon, {\rm w}^\varepsilon] \ \longrightarrow \ [\tilde{\bf u},\tilde{\rm w}] \
\mbox{ strongly in }{\bf H}_\sigma\times H^1(\Omega).
\end{equation}
Now, we observe that
\begin{eqnarray}
\|\boldsymbol{g}_1^\varepsilon\|^2_{{\bf H}^{1/2}(\Gamma_1)}
-\|\tilde{\boldsymbol{g}}_1\|^2_{{\bf H}^{1/2}(\Gamma_1)}
&=&\langle\boldsymbol{g}_1^\varepsilon-\tilde{\boldsymbol{g}}_1,\boldsymbol{g}_1^\varepsilon\rangle_{\Gamma_1}
+\langle\tilde{\boldsymbol{g}}_1, \boldsymbol{g}_1^\varepsilon-\tilde{\boldsymbol{g}}_1\rangle_{\Gamma_1},\label{nec17}\\
\|\mbox{\it g}_2^\varepsilon\|^2_{{ H}^{1/2}(\Gamma_3)}-\|\tilde{\mbox{\it g}}_2\|^2_{{ H}^{1/2}(\Gamma_3)}
&=&\langle \mbox{\it g}_2^\varepsilon-\tilde{\mbox{\it g}}_2, \mbox{\it g}_2^\varepsilon\rangle_{\Gamma_3}
+\langle\tilde{\mbox{\it g}}_2, \mbox{\it g}_2^\varepsilon-\tilde{\mbox{\it g}}_2\rangle_{\Gamma_3}.\label{nec18}
\end{eqnarray}
Taking into account that $\mathcal{U}_1\subset{\bf H}^{1/2}_{00}(\Gamma_1)\subset{\bf L}^2(\Gamma_1)$
and $\mathcal{U}_2\subset{H}^{1/2}_{00}(\Gamma_3)\subset{L}^2(\Gamma_3)$, then
from (\ref{nec14})-(\ref{nec14a}), as $\varepsilon \rightarrow 0$,  we deduce that
$$
\langle\boldsymbol{g}_1^\varepsilon-\tilde{\boldsymbol{g}}_1,\boldsymbol{g}_1^\varepsilon\rangle_{\Gamma_1}
+\langle\tilde{\boldsymbol{g}}_1, \boldsymbol{g}_1^\varepsilon-\tilde{\boldsymbol{g}}_1\rangle_{\Gamma_1}\rightarrow0,\quad
\langle \mbox{\it g}_2^\varepsilon-\tilde{\mbox{\it g}}_2, \mbox{\it g}_2^\varepsilon\rangle_{\Gamma_3}
+\langle\tilde{\mbox{\it g}}_2, \mbox{\it g}_2^\varepsilon-\tilde{\mbox{\it g}}_2\rangle_{\Gamma_3}\rightarrow0,
$$
which, together with (\ref{nec17})-(\ref{nec18}), implies that
$$
\|{\boldsymbol g}_1^\varepsilon\|^2_{{\bf H}^{1/2}(\Gamma_1)}\rightarrow
\|\tilde{\boldsymbol g}_1\|^2_{{\bf H}^{1/2}(\Gamma_1)}\
\mbox{ and }\ \|\mbox{\it g}_2^\varepsilon\|^2_{{ H}^{1/2}(\Gamma_3)}\rightarrow
\|\tilde{\mbox{\it g}}_2\|^2_{{H}^{1/2}(\Gamma_3)}.
$$
In particular,
\begin{equation}\label{nec21-1}
\limsup_{\varepsilon\rightarrow0}\|\boldsymbol{ g}_1^\varepsilon\|_{{\bf H}^{1/2}(\Gamma\setminus\Gamma_2)}=\|\tilde{\boldsymbol{ g}}_1\|_{{\bf H}^{1/2}(\Gamma\setminus\Gamma_2)}\ \mbox{ and }\
\limsup_{\varepsilon\rightarrow0}\|\mbox{\it g}^\varepsilon_2\|_{{H}^{1/2}(\Gamma)}=\|\tilde{ \mbox{\it g}}_2\|_{{H}^{1/2}(\Gamma)}.
\end{equation}

Since $\mathcal{U}_1\times\mathcal{U}_2$ is a closed set of ${\bf H}^{1/2}_{00}(\Gamma_1)\times H^{1/2}(\Gamma_3),$ from (\ref{nec21-1}) and Proposition 3.32 of \cite{brezis}
p.78,  as $\varepsilon\rightarrow0,$ we obtain that
\begin{equation}\label{strong2}
[\boldsymbol{g}^\varepsilon_1, \mbox{\it g}_2^\varepsilon] \ \longrightarrow \ [\tilde{\boldsymbol g}_1, \tilde{\mbox{\it g}}_2]
\ \mbox{ strongly in }\mathcal{U}_1\times\mathcal{U}_2.
\end{equation}
Therefore, from (\ref{strong1}) and (\ref{strong2}), we conclude (\ref{nec4}).

In order to prove (\ref{nec8}), from (\ref{nec13})-(\ref{nec14}) and, since that $J$ is weakly lower semicontinuous
on  $\mathcal{S}_{ad}$, we obtain
$$
J[\tilde{\bf u},\tilde{\rm w},\tilde{\boldsymbol{g}}_1,\tilde{\mbox{\it g}}_2]\leq\liminf_{\varepsilon\rightarrow 0} J
({\boldsymbol s}^\varepsilon)\leq\limsup_{\varepsilon\rightarrow 0} J({\boldsymbol s}^\varepsilon)
\leq J[\tilde{\bf u},\tilde{\rm w},\tilde{\boldsymbol{g}}_1,\tilde{\mbox{\it g}}_2],
$$
which implies
\begin{equation}\label{nec22}
\lim_{\varepsilon\rightarrow 0}J({\boldsymbol s}^\varepsilon)
=J[\tilde{\bf u},\tilde{\rm w},\tilde{\boldsymbol{g}}_1,\tilde{\mbox{\it g}}_2].
\end{equation}
Then, from (\ref{nec13}) and (\ref{nec22}), we obtain
$$
J[\tilde{\bf u},\tilde{\rm w},\tilde{\boldsymbol{g}}_1,\tilde{\mbox{\it g}}_2]=\lim_{\varepsilon\rightarrow 0}
J[{\bf u}^\varepsilon,{\rm w}^\varepsilon,\boldsymbol{g}_1^\varepsilon,\mbox{\it g}_2^\varepsilon]\leq \lim_{\varepsilon\rightarrow 0}
J_\varepsilon[{\bf u}^\varepsilon,{\rm w}^\varepsilon,\boldsymbol{g}_1^\varepsilon,\mbox{\it g}_2^\varepsilon]\leq
J[\tilde{\bf u},\tilde{\rm w},\tilde{\boldsymbol{g}}_1,\tilde{\mbox{\it g}}_2],
$$
which implies (\ref{nec8}).\hfill$\diamond$\vskip12pt

\subsection{Existence of Lagrange Multipliers and Adjoint Equations}
 For simplicity, we consider the following operators
\begin{equation}\label{operators}
\left\{
\begin{array}{rcl}
K:{\bf H}_\sigma\times {\bf H}_\sigma \times {\bf H}_\sigma\times \tilde{\bf H}_\sigma \rightarrow {\bf H}',&&
\tilde{K}:{\bf H}_\sigma\times {\bf H}_\sigma \times H^1(\Omega)\times H^1_0(\Omega) \rightarrow {\bf H}',\\
B^T:{\bf H}_\sigma\times {\bf H}_\sigma\times \tilde{\bf H}_\sigma \rightarrow {\bf H}', &&
E:{\bf H}_\sigma \times H^1(\Omega)\times H^1_0(\Omega) \rightarrow {\bf H}',\\
\tilde{E}:{\bf H}_\sigma \times {\bf H}_\sigma \times H^1_0(\Omega) \rightarrow (H^1(\Omega))',
&&F':{\bf H}_\sigma\times \tilde{\bf H}_\sigma \rightarrow {\bf H}', \\
G':{\bf H}_\sigma\times  H^1_0(\Omega) \rightarrow {\bf H}',&&
\end{array}
\right.
\end{equation}
defined by
\begin{equation}\label{operadores}
\left\{
\begin{array}{rcl}
\langle K[{\bf u}, {\bf u}, {\bf u}, \mbox{\boldmath$\lambda$}],{\bf v}\rangle
&=&(\eta'(N{\bf u})(N{\bf v}){\bf u}\cdot\nabla{\bf u},{\boldsymbol\lambda}) \ \ \forall {\bf v} \in {\bf H}_\sigma,\\
\langle \tilde{K}[{\bf u}, {\bf u}, {\rm w}, \kappa],{\bf v}\rangle
&=&(\eta'(N{\bf u})(N{\bf v}){\bf u}\cdot\nabla {\rm w}, \kappa) \ \ \forall {\bf v} \in {\bf H}_\sigma,\\
\langle B^T[{\bf u}, {\bf u}, {\boldsymbol\lambda}],{\bf v}\rangle
&=& (\eta(N{\bf u})({\bf u}\cdot\nabla{\bf v}
+{\bf v}\cdot\nabla{\bf u}),{\boldsymbol\lambda}) \ \ \forall {\bf v} \in {\bf H}_\sigma,\\
\langle E[{\bf u}, {\rm w}, \kappa], {\bf v}\rangle&=&(\eta(N{\bf u}){\bf v}\cdot\nabla {\rm w}, \kappa)
\ \ \forall {\bf v} \in {\bf H}_\sigma,\\
(\tilde E[{\bf u}, {\bf u}, \kappa], {\rm z})&=&(\eta(N{\bf u}){\bf u}\cdot\nabla {\rm z}, \kappa)
\ \ \forall {\rm z} \in H^1(\Omega),\\
\langle F'[{\bf u},{\boldsymbol\lambda}], {\bf v}\rangle
&=&(\eta'(N{\bf u})(N{\bf v}){\bf f},{\boldsymbol\lambda}) \ \ \forall {\bf v} \in {\bf H}_\sigma,\\
\langle G'[{\bf u}, \kappa], {\bf v}\rangle&=&(\eta'(N{\bf u})(N{\bf v}){\rm g}, \kappa)
\ \ \forall {\bf v} \in {\bf H}_\sigma,
\end{array}
\right.
\end{equation}
where $\eta'$ denotes the first derivative of $\eta$.
\begin{theorem}\label{theo_mul-1}
Let ${\bf f}\in{\bf L}^2(\Omega)$, ${\rm g}\in L^2(\Omega)$, ${\bf u}_0\in{\bf H}^{1/2}_{00}(\Gamma_0),$
$\boldsymbol{g}_1\in\mathcal{U}_1$, ${\rm w}_0\in H^{1/2}_{00}(\Gamma_0)$, $g_2\in\mathcal{U}_2$ and
$\eta\in C^1(\mathbb{R})$. Then, for any optimal solution $[{\bf u}^\varepsilon, {\rm w}^\varepsilon,
\boldsymbol{g}_1^\varepsilon, g_2^\varepsilon]\in\mathbb{X}$ of problem
(\ref{nec1})-(\ref{nec2}) there exist Lagrange multipliers  $[\boldsymbol{\lambda}^\varepsilon,\phi^\varepsilon,
\boldsymbol{\xi}^\varepsilon, \vartheta^\varepsilon]\in\tilde{\bf H}_\sigma\times H^1_0(\Omega)\times{\bf H}^{-1/2}_{00}(\Gamma\setminus\Gamma_2)\times H^{-1/2}(\Gamma)$
given by

\begin{eqnarray}
{\boldsymbol\lambda}^\varepsilon&=&\frac{1}{\varepsilon}A^{-1}(\mu_1 A{\bf u}^\varepsilon
+B[{\bf u}^\varepsilon,{\bf u}^\varepsilon,{\bf u}^\varepsilon]
-2\mu_r{\rm rot}\,{\rm w}^\varepsilon-F({\bf u}^\varepsilon))\in\tilde{\bf H}_\sigma,\label{adj1}\\
\phi^\varepsilon&=&\frac{1}{\varepsilon}\tilde{A}^{-1}(\mu_2 \tilde{A} {\rm w}^\varepsilon
+\tilde{B}[{\bf u}^\varepsilon,{\bf u}^\varepsilon,{\rm w}^\varepsilon]+4\mu_r{\rm w}^\varepsilon
-2\mu_r{\rm rot}\,{\bf u}^\varepsilon-G({\bf u}^\varepsilon))\in H^1_0(\Omega),\label{adj2}\\
{\boldsymbol\xi}^\varepsilon&=&\frac{1}{\varepsilon}({\bf u}^\varepsilon_{|_{\Gamma\setminus\Gamma_2}}-
{\bf u}_{\boldsymbol{g}^\varepsilon_1})\in {\bf L}^2(\Gamma\setminus\Gamma_2)\subset{\bf H}^{-1/2}_{00}(\Gamma\setminus\Gamma_2),\label{adj3}\\
\vartheta^\varepsilon&=&\frac{1}{\varepsilon}({\rm w}^\varepsilon_{|_\Gamma}-{\rm w}_{g^\varepsilon_2})
\in { L}^2(\Gamma)\subset H^{-1/2}(\Gamma),\label{adj4}
\end{eqnarray}
which satisfies the following system
\begin{eqnarray}
&&\mu_1\langle A{\boldsymbol\lambda}^\varepsilon,{\bf v}\rangle+\langle K[{\bf u}^\varepsilon,{\bf u}^\varepsilon,{\bf u}^\varepsilon,
{\boldsymbol\lambda}^\varepsilon],{\bf v}\rangle+\langle B^T[{\bf u}^\varepsilon,{\bf u}^\varepsilon,
{\boldsymbol\lambda}^\varepsilon],{\bf v}\rangle+\langle\tilde{K}[{\bf u}^\varepsilon,
{\bf u}^\varepsilon, {\rm w}^\varepsilon,\phi^\varepsilon],{\bf v}\rangle\nonumber\\
&&\ \ \ +\langle E[{\bf u}^\varepsilon,{\rm w}^\varepsilon,\phi^\varepsilon],{\bf v}\rangle+\langle J'_{\bf u}({\bf u}^\varepsilon),{\bf v}\rangle+
\langle {\boldsymbol\xi}^\varepsilon,{\bf v}\rangle_{\Gamma\setminus\Gamma_2}+({\bf u}^\varepsilon-\tilde{\bf u},{\bf v})_{{\bf H}_\sigma}\nonumber\\
&&\hspace{2cm}=2\mu_r({\rm rot}\,\phi^\varepsilon,{\bf v})+\langle F'[{\bf u}^\varepsilon,{\boldsymbol\lambda}^\varepsilon],{\bf v}\rangle
+\langle G'[{\bf u}^\varepsilon,\phi^\varepsilon],{\bf v}\rangle \ \ \forall {\bf v}\in {\bf H}_\sigma,\label{adj5}\\
&&\mu_2(\tilde{A}\phi^\varepsilon,{\rm z})+\langle\tilde E[{\bf u}^\varepsilon,{\bf u}^\varepsilon,\phi^\varepsilon],{\rm z}\rangle
+4\mu_r(\phi^\varepsilon,{\rm z})+({\rm w}^\varepsilon-\tilde{{\rm w}},{\rm z})_{H^1}+\beta_3( {\rm w}^\varepsilon-{\rm w}_d,{\rm z})+\langle\vartheta^\varepsilon, {\rm z}\rangle_{\Gamma}\qquad \nonumber\\
&&\hspace{2cm}=2\mu_r({\rm rot}\,{\boldsymbol\lambda}^\varepsilon,{\rm z})\ \ \forall {\rm z} \in H^1(\Omega),\label{adj6}\\
&&\beta_5\langle\boldsymbol{g}_1^\varepsilon, \boldsymbol{g}_1-\boldsymbol{g}_1^\varepsilon\rangle_{\Gamma_1}
+\langle\boldsymbol{g}_1^\varepsilon-\tilde{\boldsymbol{g}}_1, \boldsymbol{g}_1-\boldsymbol{g}_1^\varepsilon\rangle_{\Gamma_1}
-\langle{\boldsymbol\xi}^\varepsilon,\boldsymbol{g}_1-\boldsymbol{g}_1^\varepsilon
\rangle_{\Gamma_1}\geq 0 \ \ \ \forall \boldsymbol{g}_1 \in \mathcal{U}_1,\label{adj7}\\
&&\beta_6\langle g_2^\varepsilon, g_2- g_2^\varepsilon\rangle_{\Gamma_3}
+\langle g_2^\varepsilon-\tilde{g}_2, g_2- g_2^\varepsilon\rangle_{\Gamma_3}
-\langle\vartheta^\varepsilon, g_2-g_2^\varepsilon\rangle_{\Gamma_3}\geq 0 \ \ \
\forall g_2 \in \mathcal{U}_2,\label{adj8}
\end{eqnarray}
where
\begin{equation}\label{adj9}
\langle J'_{\bf u}({\bf u}^\varepsilon), {\bf v}\rangle :=\beta_1({\rm rot}\,{{\bf u}}^\varepsilon,{\rm rot}\,{\bf v})
+\beta_2({\bf u}^\varepsilon-{\bf u}_d,{\bf v})+\beta_4(\eta(N{\bf u}^\varepsilon)-\rho_d,\eta'(N{\bf u}^\varepsilon)N{\bf v}).
\end{equation}
Moreover, there exists a constant $C>0$, independent of $\varepsilon$, such that
\begin{eqnarray}
\|\boldsymbol{\xi}^\varepsilon\|_{{\bf H}^{-1/2}_{00}(\Gamma\setminus\Gamma_2)}&\leq&
C+C\|[\boldsymbol{\lambda}^\varepsilon,\phi^\varepsilon]\|_{\tilde{\bf H}_\sigma\times H^1_0},\label{adj9-1}\\
\|\vartheta^\varepsilon\|_{H^{-1/2}(\Gamma)}&\leq&C+C
\|[\boldsymbol{\lambda}^\varepsilon,\phi^\varepsilon]\|_{\tilde{\bf H}_\sigma\times H^1_0}.\label{adj9-2}
\end{eqnarray}

\end{theorem}
\textit{Proof.}
We introduce the function $\mathcal{F}:[0,1]\times[0,1]\times[0,1]\times[0,1]\longrightarrow\mathbb{R}$ defined by

\begin{equation}\label{adj10}
\mathcal{F}[\zeta_1,\zeta_2,\zeta_3,\zeta_4]=J_\varepsilon[{\bf u}^\varepsilon+\zeta_1{\bf v},{\rm w}^\varepsilon+\zeta_2{\rm z},
\boldsymbol{g}_1^\varepsilon+\zeta_3(\boldsymbol{g}_1-\boldsymbol{g}_1^\varepsilon),g_2^\varepsilon+\zeta_4(g_2-g_2^\varepsilon)],
\end{equation}
where $[{\bf v},{\rm z},\boldsymbol{g}_1,g_2]\in\mathbb{X}$.

Since the function $\mathcal{F}$ attains its minimum at ${\bf 0}=[0,0,0,0]$ and
$\mathcal{U}_1\times\mathcal{U}_2$ is convex, we have
\begin{equation}\label{adj11}
\frac{\partial \mathcal{F}}{\partial\zeta_1}({\bf 0})=0,\quad \frac{\partial \mathcal{F}}{\partial\zeta_2}({\bf 0})=0,\quad
\frac{\partial \mathcal{F}}{\partial\zeta_3}({\bf 0})\geq0,\quad \frac{\partial \mathcal{F}}{\partial\zeta_4}({\bf 0})\geq0.
\end{equation}
Therefore, from (\ref{adj11}), and the definitions of
 $\boldsymbol{\lambda}^\varepsilon,\phi^\varepsilon,
\boldsymbol{\xi}^\varepsilon, \vartheta^\varepsilon$ given in (\ref{adj1})-(\ref{adj4}),
we obtain the system (\ref{adj5})-(\ref{adj8}).

Now we will prove inequality (\ref{adj9-1}).   From (\ref{adj5}) we get
\begin{eqnarray}\label{adj14}
|\langle{\boldsymbol\xi}^\varepsilon,{\bf v}\rangle_{\Gamma\setminus\Gamma_2}|&\leq&\mu_1|\langle A{\boldsymbol\lambda}^\varepsilon,{\bf v}\rangle|
+|\langle K[{\bf u}^\varepsilon,{\bf u}^\varepsilon,{\bf u}^\varepsilon,{\boldsymbol\lambda}^\varepsilon],{\bf v}\rangle|
+|\langle B^T[{\bf u}^\varepsilon,{\bf u}^\varepsilon,{\boldsymbol\lambda}^\varepsilon],{\bf v}\rangle|\nonumber\\
&&+|\langle \tilde{K}[{\bf u}^\varepsilon,{\bf u}^\varepsilon,{\rm w}^\varepsilon,\phi^\varepsilon],{\bf v}\rangle|
+|\langle E[{\bf u}^\varepsilon,{\rm w}^\varepsilon,\phi^\varepsilon],{\bf v}\rangle |+|\langle J'_{\bf u}({\bf u}^\varepsilon),{\bf v}\rangle|\nonumber\\
&&+|({\bf u}^\varepsilon-\tilde{\bf u}, {\bf v})_{{\bf H}_\sigma}|+2\mu_r|({\rm rot}\,\phi^\varepsilon,{\bf v})|
+|\langle F'[{\bf u}^\varepsilon,{\boldsymbol\lambda}^\varepsilon],{\bf v}\rangle|\nonumber\\
&&+|\langle G'[{\bf u}^\varepsilon,\phi^\varepsilon],{\bf v}\rangle|.
\end{eqnarray}
We shall find bounds for the terms on right hand side of (\ref{adj14}). By using the H\"{o}lder inequality and
observing that $\|{\bf v}\|_{{\bf H}_\sigma}\leq C\|{\bf v}\|_{\tilde{\bf H}_\sigma}$, we obtain
\begin{eqnarray}
\mu_1|\langle A{\boldsymbol\lambda}^\varepsilon,{\bf v}\rangle |&\leq&2\mu_1|(D({\boldsymbol\lambda}^\varepsilon),D({\bf v}))|
+2\alpha\mu_1\left|\int_{\Gamma_2}\boldsymbol{\lambda}^\varepsilon\cdot{\bf v}\,d\Gamma\right|\nonumber\\
&\leq &2\mu_1\|D(\boldsymbol{\lambda}^\varepsilon)\|\|D({\bf v})\|
+2\alpha\mu_1\|\boldsymbol{\lambda}^\varepsilon\|_{{\bf L}^2(\Gamma)}\|{\bf v}\|_{{\bf L}^2(\Gamma)}\nonumber\\
&\leq&C\|{\bf v}\|_{{\bf H}_\sigma}\|{\boldsymbol\lambda}^\varepsilon\|_{\tilde{\bf H}_\sigma}
+2\mu_1\alpha C\|\boldsymbol{\lambda}^\varepsilon\|_{\tilde{\bf H}_\sigma}\|{\bf v}\|_{{\bf H}_\sigma}
\leq C\|{\bf v}\|_{{\bf H}_\sigma}\|{\boldsymbol\lambda}^\varepsilon\|_{\tilde{\bf H}_\sigma},\label{adj15}\\
2\mu_r|({\rm rot}\,\phi^\varepsilon,{\bf v})|&\leq& 2\mu_r\|{\rm rot}\,\phi^\varepsilon\|\|{\bf v}\|
\leq C\mu_r\|\phi^\varepsilon\|_{H^1_0}\|{\bf v}\|_{{\bf H}_\sigma}
\ \leq \ C\|{\bf v}\|_{{\bf H}_\sigma}\|\phi^\varepsilon\|_{H^1_0}.\label{adj15a}
\end{eqnarray}
Now, by using the H\"{o}lder inequality and the fact that $\|{\bf u}^\varepsilon\|_{{\bf H}_\sigma}\leq C$,
where $C$ is independent of $\varepsilon$, we obtain
\begin{eqnarray}\label{adj16}
|({\bf u}^\varepsilon-\tilde{\bf u}, {\bf v})_{{\bf H}_\sigma}|
&\leq&C\|{\bf u}^\varepsilon-\tilde{\bf u}\|_{{\bf H}_\sigma}\|{\bf v}\|_{{\bf H}_\sigma}
\ \leq C\|{\bf v}\|_{{\bf H}_\sigma}.
\end{eqnarray}
From (\ref{nec4}), the definition of $N$ and $\eta$ given in Subsection \ref{sub2.1}, we deduce that
\begin{equation}\label{adj17}
\|\eta(N({\bf u}^\varepsilon)\|_{\infty}\leq C,\
 \|\eta'(N({\bf u}^\varepsilon)\|_{\infty}\leq C, \
 \|N{\bf v}\|_{\infty}\leq C\|N{\bf v}\|_{H^2}\leq C\|{\bf v}\|_{{\bf H}_\sigma},
\end{equation}
where $C$ is a constant independent of $\varepsilon$.
By the H\"{o}lder inequality, (\ref{adj17}) and the definition of $\langle J'_{\bf u}({\bf u}^\varepsilon),
{\bf v}\rangle$ given in (\ref{adj9}), we obtain
\begin{eqnarray}\label{adj18}
|\langle J'_{\bf u}({\bf u}^\varepsilon),{\bf v}\rangle|&\leq&\beta_1|({\rm rot}\,{\bf u}^\varepsilon,
{\rm rot}\,{\bf v})|+\beta_2|({\bf u}^\varepsilon-{\bf u}_d,{\bf v})|+\beta_4|(\eta(N({\bf u}^\varepsilon))
-\rho_d, \eta'(N({\bf u}^\varepsilon)) N{\bf v})|\nonumber\\
&\leq&\beta_1\|{\rm rot}\,{\bf u}^\varepsilon\|\|{\rm rot}\,{\bf v}\|+\beta_2\|{\bf u}^\varepsilon
-{\bf u}_d\|\|{\bf v}\|+\beta_4\|\eta(N{\bf u}^\varepsilon)-\rho_d\|_\infty
\|\eta'(N{\bf u}^\varepsilon)\|\|N{\bf v}\|\nonumber\\
&\leq & \beta_1 C\|{\bf v}\|_{{\bf H}_\sigma}+\beta_2C \|{\bf v}\|_{{\bf H}_\sigma}
+\beta_4C \|{\bf v}\|_{{\bf H}_\sigma} \ \leq \ C\|{\bf v}\|_{{\bf H}_\sigma}.
\end{eqnarray}
Also, by the H\"{o}lder and Poincar\'e inequalities, (\ref{adj17}) and the definition of the operators $B^T$, $E$, $K$, $\tilde{K}$,
$F'$, and $G'$ given in (\ref{operadores}), we obtain

\begin{eqnarray}
|\langle B^T[{\bf u}^\varepsilon,{\bf u}^\varepsilon,{\boldsymbol\lambda}^\varepsilon],{\bf v}\rangle |&=&
|\langle B[{\bf u}^\varepsilon,{\bf u}^\varepsilon,{\bf v}]+B[{\bf u}^\varepsilon,{\bf v},{\bf u}^\varepsilon], {\boldsymbol\lambda}^\varepsilon\rangle|\nonumber\\
&\leq&|(\eta(N{\bf u}^\varepsilon){\bf u}^\varepsilon\cdot\nabla{\bf v},\boldsymbol\lambda^\varepsilon)|+
|(\eta(N{\bf u}^\varepsilon){\bf v}\cdot\nabla{\bf u}^\varepsilon,\boldsymbol\lambda^\varepsilon)|\nonumber\\
&\leq &\|\eta(N{\bf u}^\varepsilon)\|_{\infty}(\|{\bf u}^\varepsilon\|_3\|\nabla{\bf v}\|
+\|{\bf v}\|_3\|\nabla{\bf u}^\varepsilon\|)\|\boldsymbol\lambda^\varepsilon\|_6 \nonumber\\
&\leq&C\|{\bf u}^\varepsilon\|_{{\bf H}_\sigma}\|{\bf v}\|_{{\bf H}_\sigma}\|\boldsymbol\lambda^\varepsilon\|_{\tilde{\bf H}_\sigma}
\ \leq \ C\|{\bf v}\|_{{\bf H}_\sigma}\|\boldsymbol\lambda^\varepsilon\|_{\tilde{\bf H}_\sigma},\label{adj19}
\end{eqnarray}
\begin{eqnarray}
|\langle E[{\bf u}^\varepsilon,{\rm w}^\varepsilon,\phi^\varepsilon],{\bf v}\rangle|&=&
|(\eta(N{\bf u}^\varepsilon){\bf v}\cdot\nabla {\rm w}^\varepsilon,\phi^\varepsilon)|\nonumber\\
&\leq&\|\eta(N{\bf u}^\varepsilon)\|_{\infty}\|{\bf v}\|_3\|\nabla{\rm w}^\varepsilon\|\|\phi^\varepsilon\|_6\nonumber\\
&\leq&C\|{\bf v}\|_{{\bf H}_\sigma}\|{\rm w}^\varepsilon\|_{H^1}\|\phi^\varepsilon\|_{H^1_0}
 \ \leq \ C \|{\bf v}\|_{{\bf H}_\sigma}\|\phi^\varepsilon\|_{H^1_0},\label{adj20}
\end{eqnarray}
\begin{eqnarray}
|\langle K[{\bf u}^\varepsilon,{\bf u}^\varepsilon,{\bf u}^\varepsilon,{\boldsymbol\lambda}^\varepsilon],{\bf v}\rangle|&=&
|(\eta'(N{\bf u}^\varepsilon)(N{\bf v}){\bf u}^\varepsilon\cdot\nabla{\bf u}^\varepsilon,{\boldsymbol\lambda}^\varepsilon)|\nonumber\\
&\leq & \|\eta'(N{\bf u}^\varepsilon)\|_{\infty}\|N{\bf v}\|_{\infty}\|{\bf u}^\varepsilon\|_3
\|\nabla{\bf u}^\varepsilon\|\|\boldsymbol\lambda^\varepsilon\|_6\nonumber\\
&\leq&C\|{\bf v}\|_{{\bf H}_\sigma}\|{\bf u}^\varepsilon\|^2_{{\bf H}_\sigma}\|\boldsymbol\lambda^\varepsilon\|_{\tilde{\bf H}_\sigma}
\ \leq \ C\|{\bf v}\|_{{\bf H}_\sigma}\|\boldsymbol\lambda^\varepsilon\|_{\tilde{\bf H}_\sigma},\label{adj21}
\end{eqnarray}
\begin{eqnarray}
|\langle \tilde{K}[{\bf u}^\varepsilon,{\bf u}^\varepsilon,{\rm w}^\varepsilon,\phi^\varepsilon],{\bf v}\rangle|
&=&|(\eta'(N{\bf u}^\varepsilon)(N{\bf v}){\bf u}^\varepsilon\cdot\nabla {\rm w}^\varepsilon,\phi^\varepsilon)|\nonumber\\
&\leq&\|\eta'(N{\bf u}^\varepsilon)\|_{\infty}\|N{\bf v}\|_{\infty}\|{\bf u}^\varepsilon\|_3
\|\nabla{\rm w}^\varepsilon\|\|\phi^\varepsilon\|_6\nonumber\\
&\leq&C\|{\bf v}\|_{{\bf H}_\sigma}\|{\bf u}^\varepsilon\|_{{\bf H}_\sigma}\|{\rm w}^\varepsilon\|_{H^1}\|\phi^\varepsilon\|_{H^1_0}
 \ \leq \ C\|{\bf v}\|_{{\bf H}_\sigma}\|\phi^\varepsilon\|_{H^1_0},\label{adj22}
 \end{eqnarray}
\begin{eqnarray}
|\langle F'[{\bf u}^\varepsilon,{\boldsymbol\lambda}^\varepsilon],{\bf v}\rangle|&=&
|(\eta'(N{\bf u}^\varepsilon)(N{\bf v}){\bf f},{\boldsymbol\lambda}^\varepsilon)|\nonumber\\
&\leq&\|\eta'(N{\bf u}^\varepsilon)\|_{\infty}\|N{\bf v}\|_{\infty}\|{\bf f}\|
\|\boldsymbol\lambda^\varepsilon\|\ \leq \ C\|{\bf v}\|_{{\bf H}_\sigma}\|\boldsymbol\lambda^\varepsilon\|_{\tilde{\bf H}_\sigma}, \label{adj23}
\end{eqnarray}
\begin{eqnarray}
|\langle G'[{\bf u}^\varepsilon,\phi^\varepsilon],{\bf v}\rangle|&=&|(\eta'(N{\bf u}^\varepsilon)(N{\bf v}){\rm g},\phi^\varepsilon)|\nonumber\\
&\leq&\|\eta'(N{\bf u}^\varepsilon)\|_{\infty}\|N{\bf v}\|_{\infty}\|{\rm g}\|
\|\phi^\varepsilon\| \ \leq \ C\|{\bf v}\|_{{\bf H}_\sigma}\|\phi^\varepsilon\|_{H^1_0}.\label{adj24}
\end{eqnarray}
By substituting inequalities (\ref{adj15})-(\ref{adj16}) and (\ref{adj18})-(\ref{adj24}) in (\ref{adj14}), we obtain
\begin{eqnarray*}
|\langle{\boldsymbol\xi}^\varepsilon,{\bf v}\rangle_{\Gamma\setminus\Gamma_2}|
&\leq& C\|{\bf v}\|_{{\bf H}_\sigma}+ C \|{\bf v}\|_{{\bf H}_\sigma}(\|\boldsymbol\lambda^\varepsilon\|_{\tilde{\bf H}_\sigma}
+ \|\phi^\varepsilon\|_{H^1_0}),
\end{eqnarray*}
and then, we deduce that
\[
\|\boldsymbol\xi^\varepsilon\|_{{\bf H}^{-1/2}_{00}(\Gamma\setminus\Gamma_2)}\leq C+ C(\|\boldsymbol\lambda^\varepsilon\|_{\tilde{\bf H}_\sigma}
+ \|\phi^\varepsilon\|_{H^1_0})
\leq C+\sqrt{2}C\|[\boldsymbol\lambda^\varepsilon, \phi^\varepsilon]\|_{\tilde{\bf H}_\sigma\times H^1_0},
\]
which implies (\ref{adj9-1}). Analogously we can obtain (\ref{adj9-2}).\hfill$\diamond$\vskip12pt
\subsection{Optimality System}\label{sub5}

This subsection is devoted to obtain an optimality system to
problem (\ref{eq1c}). We first
show the existence of Lagrange multipliers.
\begin{theorem}\label{Lagrange}
Let ${\bf f}\in{\bf L}^2(\Omega)$, ${\rm g}\in L^2(\Omega)$, ${\bf u}_0\in{\bf H}^{1/2}_{00}(\Gamma_0),$
$\boldsymbol{g}_1\in\mathcal{U}_1$, ${\rm w}_0\in H^{1/2}_{00}(\Gamma_0)$, $g_2\in\mathcal{U}_2$ and
$\eta\in C^1(\mathbb{R})$. Then, for any
optimal solution $\tilde{\bf s}=[\tilde{\bf u},\tilde{\rm w},\tilde{\boldsymbol{g}}_1,\tilde{g}_2]\in\mathcal{S}_{ad}$
of problem (\ref{eq1c}) there exist Lagrange multipliers $[\lambda_0,\boldsymbol{\lambda}, \phi, \boldsymbol{\xi},\vartheta]\in
(\mathbb{R}^+\cup\{0\})\times\tilde{\bf H}_\sigma\times H^1_0(\Omega)\times{\bf H}^{-1/2}_{00}(\Gamma\setminus\Gamma_2)\times H^{-1/2}(\Gamma),$ not all zero, satisfying  the
following system:

\begin{eqnarray}
&&\mu_1\langle A{\boldsymbol\lambda},{\bf v}\rangle+\langle K[\tilde{\bf u}, \tilde{\bf u},\tilde{\bf u}, {\boldsymbol\lambda}],{\bf v}\rangle
+\langle B^T[\tilde{\bf u}, \tilde{\bf u}, {\boldsymbol\lambda}],{\bf v}\rangle+\langle \tilde{K}[\tilde{\bf u}, \tilde{\bf u},\tilde{{\rm w}}, \phi],{\bf v}\rangle
+\langle E[\tilde{\bf u}, \tilde{{\rm w}}, \phi], {\bf v}\rangle\nonumber\\
&&+\lambda_0\langle J'_{\bf u}(\tilde{\bf u}),{\bf v}\rangle+\langle{\boldsymbol\xi},{\bf v}\rangle_{\Gamma\setminus\Gamma_2}
-2\mu_r({\rm rot}\,\phi, {\bf v})= \langle F'[\tilde{\bf u},{\boldsymbol\lambda}],{\bf v}\rangle
+\langle G'[\tilde{\bf u},\phi],{\bf v}\rangle \ \ \forall {\bf v} \in {\bf H}_\sigma,\label{opt1-1}\\
&&\mu_2\langle\tilde{A}\phi,{\rm z}\rangle+\langle\tilde E[\tilde{\bf u}, \tilde{\bf u},\phi],{\rm z}\rangle+4\mu_r(\phi,{\rm z})
+\lambda_0\beta_3(\tilde{{\rm w}}-{\rm w}_d,{\rm z})+\langle  \vartheta,{\rm z}\rangle_{\Gamma}\nonumber\\
&&\hspace{6.7cm}=2\mu_r({\rm rot}\, {\boldsymbol\lambda},{\rm z}) \ \ \forall {\rm z} \in H^1(\Omega),\label{opt2-1}\\
&&\hspace{.7cm}\lambda_0\beta_5\langle\tilde{\boldsymbol{g}}_1, \boldsymbol{g}_1-\tilde{\boldsymbol{g}}_1
\rangle_{\Gamma_1}-\langle {\boldsymbol\xi},\boldsymbol{g}_1-\tilde{\boldsymbol{g}}_1
\rangle_{\Gamma_1}\geq 0\quad \forall \boldsymbol{g}_1\in \mathcal{U}_1,\label{opt3}\\
&&\hspace{.8cm}\lambda_0\beta_6\langle \tilde{g}_2, {\mbox{g}}_2-\tilde{g}_2\rangle_{\Gamma_3}
-\langle\vartheta, g_2-\tilde{g}_2\rangle_{\Gamma_3}\geq 0
\quad \forall g_2 \in\mathcal{U}_2.\label{opt4}
\end{eqnarray}
\end{theorem}
\textit{Proof.} From (\ref{adj9-1})-(\ref{adj9-2}) we have
\begin{equation}\label{cotaborde}
\|\boldsymbol{\xi}^\varepsilon\|_{{\bf H}^{-1/2}_{00}(\Gamma\setminus\Gamma_2)}
+\|\vartheta^\varepsilon\|_{H^{-1/2}(\Gamma)}\leq
C+C\|[\boldsymbol{\lambda}^\varepsilon,\phi^\varepsilon]\|_{\tilde{\bf H}_\sigma\times H^1_0}.
\end{equation}
For $\|[\boldsymbol{\lambda}^\varepsilon,\phi^\varepsilon]\|_{\tilde{\bf H}_\sigma\times H^1_0}$
we have the following cases
$$
\|[\boldsymbol{\lambda}^\varepsilon,\phi^\varepsilon]\|_{\tilde{\bf H}_\sigma\times H^1_0}\leq C\quad \mbox{ or }\quad
\|[\boldsymbol{\lambda}^\varepsilon,\phi^\varepsilon]\|_{\tilde{\bf H}_\sigma\times H^1_0}\longrightarrow\infty,
\mbox{ when }\varepsilon\rightarrow0.
$$
{\underline{\bf Case 1:}} $I_\varepsilon=\|[\boldsymbol{\lambda}^\varepsilon,\phi^\varepsilon]\|_{\tilde{\bf H}_\sigma\times H^1_0}\leq C$,
where $C$ is independent of $\varepsilon.$

Since the sequence $\{[\boldsymbol{\lambda}^\varepsilon, \phi^\varepsilon]\}_{\varepsilon>0}$ is bounded in
$\tilde{\bf H}_\sigma\times H^1_0(\Omega)$, there exists $[\boldsymbol{\lambda},\phi]
\in\tilde{\bf H}_\sigma\times H^1_0(\Omega)$ and a subsequence of $\{[\boldsymbol{\lambda}^\varepsilon,
\phi^\varepsilon]\}_{\varepsilon>0}$, still denoted by
$\{[\boldsymbol{\lambda}^\varepsilon, \phi^\varepsilon]\}_{\varepsilon>0}$, such that as $\varepsilon \rightarrow 0$,
\begin{equation}\label{opt5}
[\boldsymbol{\lambda}^\varepsilon, \phi^\varepsilon]\longrightarrow[\boldsymbol{\lambda},\phi]\
\mbox{ weakly in }\tilde{\bf H}_\sigma\times H^1_0(\Omega)\ \mbox{ and strongly in }\ {\bf L}^2(\Omega)\times L^2(\Omega).
\end{equation}
Then, taking into account the convergences (\ref{nec4}) and (\ref{opt5}), as  $\varepsilon\rightarrow0,$
we can obtain
\begin{equation}\label{opt6}
\left\{
\begin{array}{rcl}
\mu_1\langle A{\boldsymbol\lambda}^\varepsilon,{\bf v}\rangle
&\rightarrow &\mu_1\langle A{\boldsymbol\lambda},{\bf v}\rangle,\\
\langle K[{\bf u}^\varepsilon,{\bf u}^\varepsilon,{\bf u}^\varepsilon,{\boldsymbol\lambda}^\varepsilon],{\bf v}\rangle
&\rightarrow &\langle K[\tilde{\bf u},\tilde{\bf u},\tilde{\bf u},{\boldsymbol\lambda}],{\bf v}\rangle,\\
\langle \tilde{K}[{\bf u}^\varepsilon,{\bf u}^\varepsilon,{\rm w}^\varepsilon,\phi^\varepsilon],{\bf v}\rangle
&\rightarrow & \langle \tilde{K}[\tilde{\bf u},\tilde{\bf u},\tilde{{\rm w}},\phi],{\bf v}\rangle,\\
\langle B^T[{\bf u}^\varepsilon,{\bf u}^\varepsilon,{\boldsymbol\lambda}^\varepsilon],{\bf v}\rangle&\rightarrow &
\langle B^T[\tilde{\bf u},\tilde{\bf u},{\boldsymbol\lambda}],{\bf v}\rangle,\\
\langle E[{\bf u}^\varepsilon,{\rm w}^\varepsilon,\phi^\varepsilon],{\bf v}\rangle&\rightarrow &
\langle E[\tilde{{\bf u}},\tilde{{\rm w}},\phi],{\bf v}\rangle,\\
2\mu_r({\rm rot}\,\phi^\varepsilon, {\bf v})&\rightarrow &2\mu_r({\rm rot}\,\phi, {\bf v}),\\
\langle F'[{\bf u}^\varepsilon,{\boldsymbol\lambda}^\varepsilon],{\bf v}\rangle,
&\rightarrow &\langle F'[\tilde{\bf u},\boldsymbol\lambda],{\bf v}\rangle,\\
\langle G'[{\bf u}^\varepsilon,\phi^\varepsilon],{\bf v}\rangle&\rightarrow &\langle G'[\tilde{{\bf u}},\phi],{\bf v}\rangle,\\
\mu_2\langle\tilde{A}\phi^\varepsilon,{\rm z}\rangle \ &\rightarrow & \
\mu_2\langle\tilde{A}\phi,{\rm z}\rangle,\\
\langle\tilde E[{\bf u}^\varepsilon,{\bf u}^\varepsilon,\phi^\varepsilon],{\rm z}\rangle
\ &\rightarrow & \ \langle \tilde E[\tilde{\bf u},\tilde{\bf u},\phi],{\rm z}\rangle,\\
4\mu_r(\phi^\varepsilon,{\rm z})&\rightarrow&4\mu_r(\phi,{\rm z}),\\
\beta_3({\rm w}^\varepsilon-{\rm w}_d,{\rm z})&\rightarrow&\beta_3(\tilde{\rm w}-{\rm w}_d,{\rm z}),\\
2\mu_r({\rm rot}\,{\boldsymbol\lambda}^\varepsilon,{\rm z}) &\rightarrow & 2\mu_r({\rm rot}\,{\boldsymbol\lambda},{\rm z}),\\
\langle J'_{\bf u}({\bf u}^\varepsilon),{\bf v}\rangle &\rightarrow & \langle J'_{\bf u}(\tilde{\bf u}),{\bf v}\rangle,
\end{array}
\right.
\end{equation}
for all $[{\bf v},{\rm z}]\in{\bf H}_\sigma\times H^1(\Omega).$ In (\ref{opt6}), the operators $A, \tilde{A}, K, \tilde{K}, B^T, E,
\tilde{E}, F', G'$, and $J'_{\bf u}$ are defined in (\ref{operA}), (\ref{operadores}) and (\ref{adj9}).

Since $\|[\boldsymbol{\lambda}^\varepsilon, \phi^\varepsilon]\|_{\tilde{\bf H}_\sigma\times H^1_0}
\leq C$, from (\ref{cotaborde}), we have that the sequence $\{[\boldsymbol{\xi}^\varepsilon,\vartheta^\varepsilon]\}_{\varepsilon>0}$
is bounded in ${\bf H}^{-1/2}_{00}(\Gamma\setminus\Gamma_2)\times H^{-1/2}(\Gamma).$ Then, there exist
$[\boldsymbol{\xi},\vartheta]\in{\bf H}^{-1/2}_{00}(\Gamma\setminus\Gamma_2)\times H^{-1/2}(\Gamma)$
and a subsequence of $\{[\boldsymbol{\xi}^\varepsilon,\vartheta^\varepsilon]\}_{\varepsilon>0}$, still denoted by
$\{[\boldsymbol{\xi}^\varepsilon,\vartheta^\varepsilon]\}_{\varepsilon>0}$, such that as $\varepsilon \rightarrow 0$, we have
$$
[\boldsymbol{\xi}^\varepsilon,\vartheta^\varepsilon]\longrightarrow [\boldsymbol{\xi},\vartheta]\
\mbox{ weakly in }\ {\bf H}^{-1/2}_{00}(\Gamma\setminus\Gamma_2)\times H^{-1/2}(\Gamma),
$$
that is,
\begin{equation}\label{opt8}
\langle\boldsymbol{\xi}^\varepsilon,{\bf v}\rangle_{\Gamma\setminus\Gamma_2}\longrightarrow\langle\boldsymbol{\xi},{\bf v}\rangle_{\Gamma\setminus\Gamma_2}\ \forall{\bf v}\in{\bf H}_\sigma,\quad
\langle\vartheta^\varepsilon,{\rm z}\rangle_\Gamma\longrightarrow\langle\vartheta,{\rm z}\rangle_\Gamma\ \forall{\rm z}\in H^1(\Omega).
\end{equation}
Thus, observing (\ref{opt6})-(\ref{opt8}) and passing to the limit in (\ref{adj7})-(\ref{adj8}), as
$\varepsilon\rightarrow0$, we obtain the system (\ref{opt1-1})-(\ref{opt4}) with $\lambda_0=1.$
\vspace{0.3cm}

{\underline{\bf Case 2:}} $I_\varepsilon=\|[\boldsymbol{\lambda}^\varepsilon,\phi^\varepsilon]\|_{\tilde{\bf H}_\sigma
\times H^1_0}\rightarrow\infty$ as $\varepsilon\rightarrow0$.

By denoting
\begin{equation}\label{opt9}
\tilde{\boldsymbol{{\lambda}}}^\varepsilon=\frac{\boldsymbol{\lambda}^\varepsilon}{I_\varepsilon},\quad
\tilde{\phi}^\varepsilon=\frac{\phi^\varepsilon}{I_\varepsilon},
\end{equation}
for all $\varepsilon>0$, we have
\begin{equation}\label{opt10}
\|[\tilde{\boldsymbol{{\lambda}}}^\varepsilon,\tilde{\phi}^\varepsilon]\|_{\tilde{\bf H}_\sigma\times H^1_0}=1.
\end{equation}
Thus, the sequence $\{[\tilde{\boldsymbol{{\lambda}}}^\varepsilon,\tilde{\phi}^\varepsilon]\}_{\varepsilon>0}$
is bounded in $\tilde{\bf H}_\sigma\times H^1_0(\Omega).$ Then there exist $[\boldsymbol{\lambda},\phi]
\in\tilde{\bf H}_\sigma\times H^1_0(\Omega)$ and a subsequence of
$\{[\tilde{\boldsymbol{{\lambda}}}^\varepsilon,\tilde{\phi}^\varepsilon]\}_{\varepsilon>0}$,
still denoted by $\{[\tilde{\boldsymbol{{\lambda}}}^\varepsilon,\tilde{\phi}^\varepsilon]\}_{\varepsilon>0}$,
such that, as $\varepsilon \rightarrow 0,$ we have
\begin{equation}\label{opt11}
[\tilde{\boldsymbol{{\lambda}}}^\varepsilon,\tilde{\phi}^\varepsilon]\longrightarrow
[\boldsymbol{\lambda}, \phi]\ \mbox{ weakly in }\tilde{\bf H}_\sigma\times H^1_0(\Omega)\ \mbox{ and strongly in }\
{\bf L}^2(\Omega)\times L^2(\Omega).
\end{equation}
Moreover, by denoting
\begin{equation}\label{opt9a}
\tilde{\boldsymbol{\xi}}^\varepsilon=\frac{\boldsymbol{\xi}}{I_\varepsilon},\quad
\tilde{\vartheta}^\varepsilon=\frac{\vartheta^\varepsilon}{I_\varepsilon},
\end{equation}
from (\ref{adj9-1}) and (\ref{adj9-2}), we obtain that
$\|\tilde{\boldsymbol{\xi}}^\varepsilon\|_{{\bf H}^{-1/2}_{00}(\Gamma\setminus\Gamma_2)}
+\|\tilde{\vartheta}^\varepsilon\|_{H^{-1/2}(\Gamma)} \leq \frac{C}{I_\varepsilon}+C\leq C$,
which implies that $\{[\tilde{\boldsymbol{\xi}}^\varepsilon,\tilde{\vartheta}^\varepsilon]\}_{\varepsilon>0}$
is bounded in ${\bf H}^{-1/2}_{00}(\Gamma\setminus\Gamma_2)\times H^{-1/2}(\Gamma).$ Then, there exist
$[\boldsymbol{\xi},\vartheta]\in{\bf H}^{-1/2}_{00}(\Gamma\setminus\Gamma_2)\times H^{-1/2}(\Gamma)$
and a subsequence of $\{[\tilde{\boldsymbol{\xi}}^\varepsilon,\tilde{\vartheta}^\varepsilon]\}_{\varepsilon>0}$,
still denoted by $\{[\tilde{\boldsymbol{\xi}}^\varepsilon,\tilde{\vartheta}^\varepsilon]\}_{\varepsilon>0}$,
such that, as $\varepsilon \rightarrow 0$, we have
\begin{equation}\label{opt12}
\langle\tilde{\boldsymbol{\xi}}^\varepsilon,{\bf v}\rangle_{\Gamma\setminus\Gamma_2}\longrightarrow\langle\boldsymbol{\xi},{\bf v}\rangle_{\Gamma\setminus\Gamma_2}
\ \forall{\bf v}\in{\bf H}_\sigma,\qquad
\langle\tilde{\vartheta}^\varepsilon,{\rm z}\rangle_\Gamma\longrightarrow\langle\vartheta,{\rm z}\rangle_\Gamma\ \forall {\rm z}\in H^1(\Omega).
\end{equation}
Observing (\ref{opt9}), (\ref{opt9a}), and dividing the terms of the system (\ref{adj5})-(\ref{adj8})
by $I_\varepsilon$, we obtain
\begin{eqnarray}
&&\mu_1\langle A\tilde{\boldsymbol\lambda}^\varepsilon,{\bf v}\rangle
+\langle K[{\bf u}^\varepsilon,{\bf u}^\varepsilon,{\bf u}^\varepsilon,
\tilde{\boldsymbol\lambda}^\varepsilon],{\bf v}\rangle
+\langle B^T[{\bf u}^\varepsilon,{\bf u}^\varepsilon,\tilde{\boldsymbol\lambda}^\varepsilon],
{\bf v}\rangle+\langle\tilde{K}[{\bf u}^\varepsilon,{\bf u}^\varepsilon,{\rm w}^\varepsilon,
\tilde{\phi}^\varepsilon],{\bf v}\rangle\nonumber\\
&&\hspace{1cm}+\langle E[{\bf u}^\varepsilon,{\rm w}^\varepsilon,\tilde{\phi}^\varepsilon],{\bf v}\rangle
+\frac{1}{I_\varepsilon}\langle J'_{\bf u}({\bf u}^\varepsilon),{\bf v}\rangle
+ \langle\tilde{\boldsymbol\xi}^\varepsilon,{\bf v}\rangle_{\Gamma\setminus\Gamma_2}
+\frac{1}{I_\varepsilon}({\bf u}^\varepsilon-\tilde{\bf u},{\bf v})_{{\bf H}_\sigma}\nonumber\\
&&\hspace{1cm}=2\mu_r({\rm rot}\,\tilde{\phi}^\varepsilon, {\bf v})
+\langle F'[{\bf u}^\varepsilon,\tilde{\boldsymbol\lambda}^\varepsilon],{\bf v}\rangle
+\langle G'[{\bf u}^\varepsilon,\tilde{\phi}^\varepsilon],{\bf v}\rangle,\label{opt13}\\
&&\mu_2\langle\tilde{A}\tilde\phi^\varepsilon,{\rm z}\rangle
+\langle\tilde E[{\bf u}^\varepsilon,{\bf u}^\varepsilon,\tilde\phi^\varepsilon],{\rm z}\rangle
+4\mu_r(\tilde\phi^\varepsilon,{\rm z})+\frac{1}{I_\varepsilon}({\rm w}^\varepsilon-\tilde{{\rm w}},{\rm z})_{H^1}
+\beta_3\frac{1}{I_\varepsilon}( {\rm w}^\varepsilon-{\rm w}_d,{\rm z})\nonumber\\
&&\hspace{1cm}+\langle\tilde\vartheta^\varepsilon, {\rm z}\rangle_{\Gamma}
=2\mu_r({\rm rot}\,\tilde{\boldsymbol\lambda}^\varepsilon,{\rm z}),\label{opt14}\\
&&\frac{\beta_5}{I_\varepsilon}\langle{\boldsymbol g}_1^\varepsilon,
{\boldsymbol g}_1-{\boldsymbol g}_1^\varepsilon\rangle_{\Gamma_1}
+\frac{1}{I_\varepsilon}\langle{\boldsymbol g}_1^\varepsilon-\tilde{\boldsymbol g}_1, {\boldsymbol g}_1-\boldsymbol{g}_1^\varepsilon\rangle_{\Gamma_1}
-\langle\tilde{\boldsymbol\xi}^\varepsilon,{\boldsymbol g}_1-{\boldsymbol g}_1^\varepsilon
\rangle_{\Gamma_1}\geq 0,\label{opt15}\\
&&\frac{\beta_6}{I_\varepsilon}\langle g_2^\varepsilon, g_2- g_2^\varepsilon\rangle_{\Gamma_3}
+\frac{1}{I_\varepsilon}\langle g_2^\varepsilon-\tilde{g}_2, g_2- g_2^\varepsilon\rangle_{\Gamma_3}
-\langle\tilde{\vartheta}^\varepsilon, g_2-g_2^\varepsilon\rangle_{\Gamma_3}\geq 0,\label{opt16}
\end{eqnarray}
for all $[{\bf v}, {\rm z}, {\boldsymbol g}_1, g_2]  \in {\bf H}_\sigma\times H^1(\Omega)\times\mathcal{U}_1\times\mathcal{U}_2$.

Therefore, observing the convergences (\ref{opt6}), (\ref{opt11}), (\ref{opt12}), and passing to the
limit in (\ref{opt13})-(\ref{opt16}) when $\varepsilon\rightarrow0$, we obtain the system
(\ref{opt1-1})-(\ref{opt4}) with $\lambda_0=0$.

Now we only need to verify that $[0,\boldsymbol{\lambda},\phi,\boldsymbol{\xi},\vartheta]\neq{\bf 0}$.

Observing that $\langle A\tilde{\boldsymbol\lambda}^\varepsilon,\tilde{\boldsymbol\lambda}^\varepsilon\rangle
=2\|\tilde{\boldsymbol\lambda}^\varepsilon\|_{\tilde{\bf H}_\sigma}^2
+2\alpha\|\tilde{\boldsymbol{\lambda}}^\varepsilon\|^2_{{\bf L}^2(\Gamma_2)}$ and $\langle\tilde{A}\tilde\phi^\varepsilon,\tilde\phi^\varepsilon\rangle
=\|\tilde\phi^\varepsilon\|_{{ H}^1_0}^2$, and replacing ${\bf v}=\tilde{\boldsymbol{\lambda}}^\varepsilon$ in (\ref{opt13}) and
${\rm z}=\tilde{\phi}^\varepsilon$ in (\ref{opt14}), we obtain
\begin{eqnarray}
2\mu_1(\|\tilde{\boldsymbol\lambda}^\varepsilon\|_{\tilde{\bf H}_\sigma}^2
+\alpha\|\tilde{\boldsymbol{\lambda}}^\varepsilon\|^2_{{\bf L}^2(\Gamma_2)})
&=&-\langle K[{\bf u}^\varepsilon,{\bf u}^\varepsilon,{\bf u}^\varepsilon,
\tilde{\boldsymbol\lambda}^\varepsilon],\tilde{\boldsymbol\lambda}^\varepsilon\rangle
-\langle B^T[{\bf u}^\varepsilon,{\bf u}^\varepsilon,\tilde{\boldsymbol\lambda}^\varepsilon],
\tilde{\boldsymbol\lambda}^\varepsilon\rangle\nonumber\\
&&-\langle\tilde{K}[{\bf u}^\varepsilon,{\bf u}^\varepsilon,{\rm w}^\varepsilon,
\tilde{\phi}^\varepsilon],\tilde{\boldsymbol\lambda}^\varepsilon\rangle
-\langle E[{\bf u}^\varepsilon,{\rm w}^\varepsilon,\tilde{\phi}^\varepsilon],
\tilde{\boldsymbol\lambda}^\varepsilon\rangle\nonumber\\
&&-\frac{1}{I_\varepsilon}\langle J'_{\bf u}({\bf u}^\varepsilon),\tilde{\boldsymbol\lambda}^\varepsilon\rangle
- \langle\tilde{\boldsymbol\xi}^\varepsilon, \tilde{\boldsymbol\lambda}^\varepsilon\rangle_{\Gamma\setminus\Gamma_2}
-\frac{1}{I_\varepsilon}({\bf u}^\varepsilon-\tilde{\bf u},\tilde{\boldsymbol\lambda}^\varepsilon)_{{\bf H}_\sigma}\nonumber\\
&&+2\mu_r({\rm rot}\,\tilde{\phi}^\varepsilon, \tilde{\boldsymbol\lambda}^\varepsilon)
+\langle F'[{\bf u}^\varepsilon,\tilde{\boldsymbol\lambda}^\varepsilon],\tilde{\boldsymbol\lambda}^\varepsilon\rangle
+\langle G'[{\bf u}^\varepsilon,\tilde{\phi}^\varepsilon],\tilde{\boldsymbol\lambda}^\varepsilon\rangle,\label{opt17}\\
\mu_2\|\tilde\phi^\varepsilon\|_{{H}^1_0}^2+4\mu_r\|\tilde{\phi}^\varepsilon\|^2&=&-\langle \tilde E[{\bf u}^\varepsilon,{\bf u}^\varepsilon,
\tilde\phi^\varepsilon],\tilde\phi^\varepsilon\rangle
-\frac{1}{I_\varepsilon}({\rm w}^\varepsilon-\tilde{{\rm w}},\tilde\phi^\varepsilon)_{H^1}\nonumber\\
&&-\beta_3\frac{1}{I_\varepsilon}( {\rm w}^\varepsilon-{\rm w}_d,\tilde\phi^\varepsilon)
-\langle\tilde\vartheta^\varepsilon, \tilde\phi^\varepsilon\rangle_{\Gamma}
+2\mu_r({\rm rot}\,\tilde{\boldsymbol\lambda}^\varepsilon,\tilde\phi^\varepsilon).\label{opt18}
\end{eqnarray}
If $[\boldsymbol{\lambda},\phi]=[{\bf 0},0]$, considering (\ref{opt5})-(\ref{opt6})
and passing to the limit in (\ref{opt17})-(\ref{opt18}) as $\varepsilon\rightarrow0$, we obtain that
$\mu_1\|\tilde{\boldsymbol\lambda}^\varepsilon\|_{\tilde{\bf H}_\sigma}^2 \longrightarrow 0$ and
$\mu_2\|\tilde\phi^\varepsilon\|_{{ H}^1_0}^2\longrightarrow 0;$ then it follows
$$
\|[\tilde{\boldsymbol\lambda}^\varepsilon, \tilde\phi^\varepsilon]\|^2_{\tilde{\bf H}_\sigma\times H^1_0}
=\|\tilde{\boldsymbol\lambda}^\varepsilon\|_{\tilde{\bf H}_\sigma}^2+ \|\tilde\phi^\varepsilon\|_{ H^1_0}^2\rightarrow 0,
$$
which contradicts the equality given in (\ref{opt10}). Therefore, we conclude that $[\boldsymbol{\lambda},\phi]\neq[{\bf 0},0]$
and consequently $[0,\boldsymbol{\lambda},\phi,\boldsymbol{\xi},\vartheta]\neq{\bf 0}$. Thus, the proof of
the theorem is finished.\hfill$\diamond$\vskip12pt

\begin{remark}\label{opt1}
From (\ref{eq8}) and (\ref{opt1-1})-(\ref{opt4}),  we obtain
the  following optimality system for problem (\ref{eq1c}) constituted by the state equations (\ref{eq8}), the adjoint equations (\ref{opt1-1})-(\ref{opt2-1}) and the optimality conditions (\ref{opt3})-(\ref{opt4}).
\vspace{0.3cm}
\end{remark}

\begin{corol}
Suppose that the assumptions of Theorem \ref{Lagrange} are satisfied and let $[\tilde{\bf u},\tilde{\rm w},\tilde{\boldsymbol{g}}_1,\tilde{g}_2]\in\mathcal{S}_{ad}$ an
optimal solution to problem (\ref{eq1c}). Let $\mu_1$ and $\mu_2$ large enough such that
\begin{equation}\label{corol1}
\tilde{\delta}> C(\|[\tilde{\bf u},\tilde{\rm w}]\|_{{\bf H}_\sigma\times H^1}+\|[\tilde{\bf u},\tilde{\rm w}]\|^2_{{\bf H}_\sigma\times H^1}+
\|{\bf f}\|+\|{\rm g}\|+\mu_r),
\end{equation}
where $\tilde{\delta}=\min\{2\mu_1,\mu_2\}$ and $C$ is a positive constant depending only on $\Omega$, $\|\eta(N\tilde{\bf u})\|_{\infty}$, and $\|\eta'(N\tilde{\bf u})\|_{\infty}$.
Then, there exists a unique $[\lambda_0,\boldsymbol{\lambda}, \phi, \boldsymbol{\xi},\vartheta]
\in(\mathbb{R}^+\cup \{0\})\times\tilde{\bf H}_\sigma\times H^1_0(\Omega)\times{\bf H}_{00}^{-1/2}(\Gamma\setminus\Gamma_2)
\times H^{-1/2}(\Gamma)$ satisfying (\ref{opt1-1})-(\ref{opt4}), with $\lambda_0=1.$
\end{corol}
\textit{Proof.} We assume that $\lambda_0=0$. Then, by setting ${\bf v}=\boldsymbol{\lambda}$ in (\ref{opt1-1}) and  ${\rm z}=\phi$ in (\ref{opt2-1}), taking into
account that $\langle{\boldsymbol\xi},{\boldsymbol\lambda}\rangle_{\Gamma\setminus\Gamma_2}=0$ and  $\langle\vartheta,\phi\rangle_{\Gamma}=0$, observing the estimates given in (\ref{adj15a}),
(\ref{adj19})-(\ref{adj24}), and definition of $A$, $\tilde{A}$, we obtain
\begin{eqnarray*}
2\mu_1\|{\boldsymbol\lambda}\|^2&\leq &C\mu_r \|[\boldsymbol{\lambda},\phi]\|^2_{\tilde{\bf H}_\sigma\times H^1_0}
+ C \|\eta(N\tilde{\bf u})\|_{\infty}(\|\tilde{\bf u}\|_{{\bf H}_\sigma}
+\|\tilde{\rm w}\|_{H^1}+\|\tilde{\bf u}\|^2_{{\bf H}_\sigma})\|[\boldsymbol{\lambda},\phi]\|^2_{\tilde{\bf H}_\sigma\times H^1_0}\\
&&+C \|\eta'(N\tilde{\bf u})\|_{\infty}(\|\tilde{\bf u}\|_{{\bf H}_\sigma}\|\tilde{\rm w}\|_{H^1}+\|{\bf f}\|+\|{\rm g}\|)
\|[\boldsymbol{\lambda},\phi]\|^2_{\tilde{\bf H}_\sigma\times H^1_0},\\
\mu_2\|\phi\|^2&\leq &C\mu_r \|[\boldsymbol{\lambda},\phi]\|^2_{\tilde{\bf H}_\sigma\times H^1_0}.
\end{eqnarray*}
Then, by adding the above inequalities, we deduce
\begin{equation*}
\tilde{\delta}\|[\boldsymbol{\lambda},\phi]\|^2_{\tilde{\bf H}_\sigma\times H^1_0}
\leq  C(\|[\tilde{\bf u},\tilde{\rm w}]\|_{{\bf H}_\sigma\times H^1}+\|[\tilde{\bf u},\tilde{\rm w}]\|^2_{{\bf H}_\sigma\times H^1}+
\|{\bf f}\|+\|{\rm g}\|+\mu_r)\|[\boldsymbol{\lambda},\phi]\|^2_{\tilde{\bf H}_\sigma\times H^1_0},
\end{equation*}
which, by applying condition (\ref{corol1}) implies that $\|[\boldsymbol{\lambda},\phi]\|^2_{\tilde{\bf H}_\sigma\times H^1_0}=0$, that is,  $\boldsymbol{\lambda}={\bf 0}$ and $\phi=0$.
In this case, the equations (\ref{opt1-1}) and (\ref{opt2-1}) can be rewritten as
$\langle{\boldsymbol\xi},{\bf v}\rangle_{\Gamma\setminus\Gamma_2}=0$ for any ${\bf v}\in {\bf H}_\sigma$ and
$\langle\vartheta,{\rm z}\rangle_\Gamma=0$ for any ${\rm z}\in H^1(\Omega)$, respectively. Hence  we  have
$\boldsymbol{\xi}={\bf 0}$ and $\vartheta=0$, which contradicts Theorem \ref{Lagrange}.
\hfill$\diamond$\vskip12pt

\begin{remark}\label{opt2}
If the Lagrange multiplier  $\lambda_0=1$, then the optimality conditions are
equivalent to
\begin{eqnarray}\label{optcond}
\langle\beta_5\tilde{\boldsymbol g}_1-{\boldsymbol\xi}, \boldsymbol{g}_1-\tilde{\boldsymbol g}_1
\rangle_{\Gamma_1}\geq 0\quad\mbox{ and }\quad
\langle \beta_6\tilde{g}_2-\vartheta, g_2-\tilde{g}_2\rangle_{\Gamma_3}\geq 0.
\end{eqnarray}
Since the set of controls $\mathcal{U}_1\times\mathcal{U}_2$ is convex, from inequalities (\ref{optcond}) we obtain
$$
\tilde{\boldsymbol{g}}_1=\mathop{\rm Proj}\limits_{\mathcal{U}_1}\left(\frac{\boldsymbol{\xi}}{\beta_5}\right)\quad
\mbox{ on }\Gamma_1,\qquad
\tilde{g}_2=\mathop{\rm  Proj}\limits_{\mathcal{U}_2}\left(\frac{\vartheta}{\beta_6}\right)\quad
\mbox{ on }\Gamma_3.
$$
\end{remark}

{\bf Acknowledgments:} The first author was supported by proyecto UTA-Mayor, 4738-17. The third author was supported by Fondo Nacional de
Financiamiento para la Ciencia, la Tecnolog\'ia y la Innovaci\'on Francisco Jos\'e de
Caldas, contrato Colciencias  FP 44842-157-2016.

\end{document}